\documentclass[12pt]{article}
\usepackage[utf8]{inputenc}
\usepackage[english]{babel}

\usepackage{algorithm}
\usepackage{algpseudocode} 
\usepackage[numbers]{natbib}
\usepackage{amsmath, amsfonts, amssymb, amsbsy, bigstrut, graphicx, enumerate, upref, longtable, comment, csquotes, booktabs, array, verbatim}
\usepackage{setspace}
\usepackage{amsthm}
\usepackage{tikz}
\usetikzlibrary{trees}
\usetikzlibrary{calc, positioning, decorations.pathmorphing, shapes.geometric}
\usepackage{float}
\usepackage{subcaption}
\usepackage{pstricks}
\usepackage{placeins}
\usepackage{cleveref}
\usepackage{enumitem}
\newcommand{\blind}{0}

\newtheorem{thm}{Theorem}[section]

\newtheorem{lmm}[thm]{Lemma}
\newtheorem{defn}[thm]{Definition}

\newtheorem{cor}[thm]{Corollary}

\theoremstyle{definition}
\newtheorem{ex}[thm]{Setting}

\DeclareMathOperator{\sign}{sign}

\newcommand{\MB}{\textbf{MB}}

\newcommand{\bbx}{\mathbf{X}}
\newcommand{\bby}{\mathbf{Y}}
\newcommand{\bbz}{\mathbf{Z}}
\newcommand{\bbw}{\mathbf{W}}
\newcommand{\bbv}{\mathbf{V}}

\newcommand{\ee}{\mathbb{E}}
\newcommand{\pp}{\mathbb{P}}
\newcommand{\rr}{\mathbb{R}}

\let\oldref\ref
\renewcommand{\ref}[1]{(\oldref{#1})}

\newcommand{\var}{\mathrm{Var}}
\newcommand{\argmax}{\operatorname{argmax}}

\begin{document}

\def\spacingset#1{\renewcommand{\baselinestretch}%
{#1}\small\normalsize} \spacingset{1}


\if0\blind
{
  \title{\bf A Fast Non-parametric Approach for Local Causal Structure Learning}
  \author{Mona Azadkia, Armeen Taeb, Peter B\"{u}hlmann \\
    {\small Department of Mathematics, ETH Z\"{u}rich}}
  \maketitle
} \fi

\bigskip
\begin{abstract}
We study the problem of causal structure learning with essentially no assumptions on the functional relationships and noise. We develop DAG-FOCI, a computationally fast algorithm for this setting that is based on the FOCI variable selection algorithm in~\cite{azadkia2021simple}. DAG-FOCI outputs the set of parents of a response variable of interest. We provide theoretical guarantees of our procedure when the underlying graph does not contain any (undirected) cycle containing the response variable of interest. Furthermore, in the absence of this assumption, we give a conservative guarantee against false positive causal claims when the set of parents is identifiable. We demonstrate the applicability of DAG-FOCI on simulated as well as a real dataset from computational biology~\cite{sachs2005causal}.
\end{abstract}

\noindent%


\section{Introduction}
\label{sec:Introduction}
Causal reasoning is important in many areas, including the sciences, decision making, and public policy. The gold standard method for determining causal relations uses randomized control experiments, but the cost and ethical concerns often hinder their feasibility. Hence, it is worth estimating causal relations from observational data, that is, data obtained from observing a system without any interventions. We study this problem in a non-parametric setting where we assume that all relevant variables are observed, excluding latent variables.  

This paper proposes a computationally efficient algorithm for local causal structure learning around a target variable of interest without restricting the functional forms, error distributions, or relying on conditional independence testing. The price to be paid for this generality is in terms of conservativeness. We provide finite sample guarantees that our algorithm will not wrongly claim a causal edge when the target variable has multiple parents: it may not lead to many causal findings, but when the underlying causal graph is a polytree, the parental set is correctly identified. From a practical perspective, in an observational study without latent variables, this is so far the only causal inference algorithm that does not rely on conditional independence testing and protects against false causal claims when there are multiple causes, in very general non-parametric models.

To develop our methodology, we consider a setting where the variables satisfy the Causal Markov property with respect to a directed acyclic graph (DAG) and are generated according to a structural equation model (SEM). As announced above, no assumptions or restrictions are imposed on the functional forms or the error distributions. We focus on recovering the set of parents of a given target variable (maybe multiple since the Markov equivalence class can entail more than one set of parents). Estimating the entire DAG can be viewed as combining the local neighborhood structures of all the nodes \cite{pellet2008using, yang2021model}. 

Our proposed algorithm, dubbed DAG-FOCI, consists of two stages. In the first stage, we identify the local structure in terms of the Markov boundary of the target variable. This is done using a greedy algorithm named FOCI, introduced in \cite{azadkia2021simple}. FOCI does not rely on any statistical testing or distributional assumptions. Furthermore, it does not require the specification of any tuning parameters and is computationally efficient. In the second stage of DAG-FOCI, we combine some pairwise independence tests with the output of the first stage (for several targets) to determine the multi-set of parents of the target variable. We emphasize here that the second stage requires only marginal statistical testing, a much easier and better-posed task than conditional independence testing.  
Our theoretical guarantees for DAG-FOCI require a condition dubbed `$\delta$-Markov Gap', which states that the target variable has larger dependencies with members of its Markov boundary than with variables outside of its Markov boundary. The $\delta$-Markov Gap condition is crucial for ensuring that the greedy algorithm in the first stage of DAG-FOCI is successful. Under such an assumption, we establish a finite sample guarantee that no false-positive causal edges are claimed by DAG-FOCI when the target variable of interest has more than one parent. Furthermore, again for finite sample cases and when the target variable has multiple parents, we show that DAG-FOCI identifies the parental set if the underlying DAG is a polytree. 

The outline of this paper is as follows. In Section~\oldref{sec:DAG-FOCI}, we present our algorithm DAG-FOCI. We explain our main theoretical results in Section~\oldref{sec:Theory}. In Section~\oldref{sec:interventional}, we describe how to incorporate interventional data into DAG-FOCI for additional identifiability. Applications to simulated and real data are presented in Section~\oldref{sec:Experiments}. We conclude with future directions in Section~\oldref{sec:conclusion}.

\subsection{Related Work}
Prior work on causal structure learning can be broadly categorized into two groups: score-based algorithms and constraint-based approaches using conditional independence testing. The latter focuses mostly on Gaussian distributions \cite{kalisch2007estimating, maathuis2009estimating, colombo2014order} and transformable to Gaussian distributions  \cite{harris2013pc};  one can also strengthen these approaches by using a kernel-based tests for conditional independence \cite{fukumizu2007kernel}. Score-based methods have been worked out for the Gaussian case \cite{chickering2002optimal, van2013ell_, Loh2014}, and there is a substantial amount of work when putting additional restrictions on the functional forms (e.g. nonlinear but additive errors \cite{hoyer2008nonlinear, peters2014causal, buhlmann2014cam, nowzohour2016score}), or on the error distributions (e.g. non-Gaussian and additive error~\cite{shimizu2006linear} or equal error noise variances \cite{peters2014identifiability}). Some methods also combine score-based and constraint-based approaches \cite{tsamardinos2006max, nandy2018high}. 

As highlighted above, in contrast to DAG-FOCI, both constraint-based and score-based methods rely on techniques or assumptions which are challenging or restrictive in the purely non-parametric case. In particular, constraint-based methods such as the PC and IC algorithms~\cite{spirtes1991algorithm, pearl1995theory} require conditional independence testing, which is known to be a hard problem in non-parametric settings~\cite{shah2020hardness}. Although there has been a large body of research exploring different ideas for conditional independence testing, there is no model-free and non-parametric method with theoretical guarantees for this purpose: one has to restrict the model class to obtain power against arbitrary alternatives. Hence, many constraint-based methods are vulnerable in applications where the data may not satisfy the properties needed for the validity of the conditional independence test. Similarly, non-parametric score-based methods are not tractable in general. In addition to requiring distributional assumptions, score-based methods can become computationally infeasible in high-dimensional settings. Improving the speed of score-based methods involves either using ad-hoc optimization ~\cite{brenner2013sparsityboost, zheng2018dags} or requiring tree structure or additive noise model~\cite{jakobsen2021structure}.

The authors in ~\cite{Gao2020APA,gao2021efficient} propose polynomial-time algorithms for structural learning in non-parametric settings, and provide finite-sample theoretical guarantees for their methods. However, their methodology and theoretical analysis require more stringent assumptions than the ones imposed in this paper. Specifically, \cite{Gao2020APA} places an `equal variance' type condition on the error distributions and \cite{gao2021efficient} requires non-source nodes to have higher conditional entropy compared to at least one of their ancestors.\footnote{In general, this assumption is not satisfied: consider the simple example $X\rightarrow Y\leftarrow Z$ and $Y = f(X,Z,\varepsilon)$ where $X, Z, \varepsilon\overset{i.i.d}{\sim}\text{Unif}[0, 1]$ and $f:[0, 1]^3\rightarrow [0, 1]$ is a measurable function.} Furthermore, neither method provides false discovery control guarantees when their assumptions are violated. 

\subsection{Notations and Definitions}\label{sec:not}
We represent a DAG $D = (V,E)$ with a node set $V = [p+1]$ and a collection of directed edges $E$. We denote a directed edge by $(i, j)$ or $i\rightarrow j$. In this case, $i$ is a \textbf{parent} of $j$, and $j$ is a \textbf{child} of $i$. If there is a directed path $i\rightarrow\cdots\rightarrow j$, then $j$ is a \textbf{descendent} of $i$, and $i$ is an \textbf{ancestor} of $j$. The \textbf{skeleton} of $D$ is the undirected graph obtained from $D$ by substituting directed edges with undirected edges. Nodes $i$, $j$, and $k$ form a \textbf{v-structure} in a DAG if $i\rightarrow k \leftarrow j$ and there is no edge between $i$ and $j$. In this case, $i$ and $j$ are \textbf{spouses} and node $k$ is called an \textbf{unshielded collider}. More generally, if $i$ and $j$ can have an edge connecting them, node $k$ is called a \textbf{collider}.

For the joint probability distribution $\mathbb{P}$ on $(\bbx_1,\cdots,\bbx_p, Y := \bbx_{p+1})\in\rr^{p+1}$, the \textbf{Causal Markov property} holds if every variable is independent of its non-descendants conditional on its parents. The \textbf{Markov equivalence class} of $D$ is the set of all DAGs that encode the same set of conditional independencies. Two DAGs are Markov equivalent if and only if they have the same skeleton and the same v-structures~\cite{verma1991equivalence}. Let $\{D_1, \cdots, D_\kappa\}$ be the Markov equivalence class of DAG $D$. We denote the set of parents, children, and spouses of $Y$ in DAG $D_i$ respectively by $\texttt{pa}_i(Y)$, $\texttt{ch}_i(Y)$, and $\texttt{sp}_i(Y)$. For simplicity in the notation and and for a special DAG of interest $D$, we drop the subscripts, e.g, $\texttt{pa}(Y)$ etc. We represent the collection of all parental sets in the Markov equivalence class by $\mathcal{P}(Y) = \{\texttt{pa}_1(Y), \cdots, \texttt{pa}_K(Y)\}$. Note that $K$ may be less than $\kappa$ as DAGs $D_k$ and $D_{l}$ in the Markov equivalence class may encode the same set of parents for $Y$, i.e. $\texttt{pa}_k(Y) = \texttt{pa}_l(Y)$. Note that if $Y$ is a source node in $D_k$, then $\texttt{pa}_k(Y)=\emptyset$. From the distribution of observational data, one can in general only identify the Markov equivalence class; in particular, this is the case in the fully non-parametric setting which we consider in this work. 

For any $i\in[p+1]$, $\bbx_{\setminus i}$ is the random vector consisting of all $\bbx_j$'s excluding $\bbx_i$. For an index set $S\subseteq[p+1]$, the random vector of corresponding features is denoted by $\bbx_S$. For $S\subseteq [p+1]\setminus\{i\}$, $S$ is a \textbf{Markov blanket} or a \textbf{sufficient set} for $\bbx_i$ if and only if $\bbx_i\perp\bbx_{[p+1]\setminus (S\cup\{i\})}\mid \bbx_S$. A minimal Markov blanket is called a \textbf{Markov boundary}, i.e., a Markov blanket is a Markov boundary if no strict subset of it is a Markov blanket. In this paper, we always assume that the Markov boundary is unique; indeed this is the case under mild assumptions \cite[Th.4.3]{pearl2014probabilistic}. We denote the unique Markov boundary of $\bbx_i$ by $\MB(\bbx_i)$. It is well-known that for any $D_i$, $\MB(Y) = \texttt{pa}_i(Y)\cup\texttt{ch}_i(Y)\cup\texttt{sp}_i(Y)$.

Finally, our methodological development and theoretical analysis is motivated by the graphical notions of tree-neighborhood and polytree, which are defined below:
\begin{defn}\label{tnb}
Node $Y$ satisfies a \textbf{tree-neighborhood} property if any cycle in the skeleton of $D$ contains at most one member of $\MB(Y)$. 
\end{defn}

\begin{defn}\label{polytree}
A \textbf{polytree} is a DAG whose skeleton is a tree. 
\end{defn}

\section{Our algorithm DAG-FOCI}
\label{sec:DAG-FOCI}

We consider the random vector $\bbx = (\bbx_1, \cdots, \bbx_p, \bbx_{p+1})\in\rr^{p+1}$ and denote by $Y = \bbx_{p+1}\in\rr$ a response variable of special interest. We posit that the underlying distribution of $\bbx$ satisfies the Causal Markov property with respect to a DAG $D^\star$ and is parameterized by the following structural equation model (SEM):
\begin{eqnarray}\label{SEM}
\bbx_i = f_i(\bbx_{i_1}, \cdots, \bbx_{i_k}, \varepsilon_i) ,\ i\in[p+1], \ \texttt{pa}(\bbx_i) = \{i_1, \cdots, i_k\} .
\end{eqnarray}
Here, $f_i: \rr^{|\texttt{pa}(\bbx_i)|} \times \rr \to \rr$ is a measurable function, and $\varepsilon_1, \ldots, \varepsilon_{p+1}$ are jointly independent noise variables from non-degenerate probability laws.

We always assume that the observed data consists of $n$ independent and identically distributed realizations from the SEM in \ref{SEM}.
Extensions to interventional data are presented in Section~\oldref{sec:interventional}.
 
As described in the introduction, our algorithm DAG-FOCI consists of two stages. In the first stage, we deploy the non-parametric variable selection procedure FOCI \cite{azadkia2021simple} to obtain the local Markov structure around the node $Y$: the Markov boundary of $Y$, denoted as $\MB(Y)$, as well as the Markov blanket of each of the variables in $\MB(Y)$. In the second stage, we combine the local Markov structures around $Y$ to identify the parents of $Y$. The description of DAG-FOCI is presented in Algorithm~\oldref{alg:overall}. 

\begin{algorithm}[h]
\setstretch{1}
\caption{DAG-FOCI}
\begin{algorithmic}[1]
\Require $n$ samples of $(\bbx_1, \cdots, \bbx_p, Y = \bbx_{p+1})$ 
\Ensure  $\hat{\mathcal{P}}_n(Y)$, the set of all possible parental sets of $Y$
\State {\bf Markov boundary search}: using FOCI (see Section \oldref{sec:FOCI}), estimate first $\widehat{\MB}(Y)$ and then also $\widehat{\MB}(\bbx_j)$ for all $\bbx_j \in \widehat{\MB}(Y)$
\State {\bf Identifying clusters of nodes}: find maximal sets $S \subseteq \widehat{\MB}(Y)$ such that for every $i,j \in S$: 
\[
\left\{ \begin{array}{l}
         i \in  \widehat{\MB}(\bbx_j);\\
         j \in  \widehat{\MB}(\bbx_i); \\
         \bbx_i \perp \bbx_j.
         \end{array} \right.
\]
Let $\hat{\mathcal{S}}_n(Y)$ be the collection of all such $S$ plus $\{\emptyset\}$. 
\State {\bf Output}: Output $\hat{\mathcal{P}}_n(Y)$ is defined in the following way:
\[
\left\{ \begin{array}{l}
         \text{if there is only one $S\in\hat{\mathcal{S}}_n(Y)$ such that $|S| > 1$, let $\hat{\mathcal{P}}_n(Y) = \{S\}$};\\
        
         \text{if for all non-empty $S\in\hat{\mathcal{S}}_n(Y)$, we have $|S| = 1$, let $\hat{\mathcal{P}}_n(Y) = \hat{\mathcal{S}}_n(Y)$};\\
         
         \text{otherwise return ``DAG-FOCI is not able to detect the parents".}
         \end{array} \right.
\]
\end{algorithmic}
\label{alg:overall}
\end{algorithm}
In Sections~\oldref{sec:FOCI} and~\oldref{Parents}, we explain Algorithm~\oldref{alg:overall} in more details. In Section~\oldref{sec:FOCI} we explain how the Markov boundary search is done using FOCI and in Section~\oldref{Parents} we elaborate on steps 2-3.

\subsection{Stage I of DAG-FOCI: Finding the Markov Boundary}\label{sec:FOCI}
For Markov boundary estimation in step 1 of Algorithm \oldref{alg:overall}, we use the model-free variable selection algorithm FOCI~\cite{azadkia2021simple}. FOCI is based on the non-parametric and distribution-free measure of dependence CODEC~\cite{azadkia2021simple}. For a random variable $Y\in\rr$ and random vectors $\bbx\in\rr^q$ and $\bbz\in\rr^s$ where $s\geq 0$ and $q\geq 1$, CODEC quantifies the dependency of $Y$ on $\bbz$ conditional on $\bbx$ as 
\begin{eqnarray}\label{TOriginal}
T(Y, \bbz\mid\bbx) = \frac{\int\ee(\var(\pp(Y\geq t\mid \bbz, \bbx)\mid\bbx))d\mu(t)}{\int\ee(\var(1_{\{Y\geq t\}}\mid\bbx))d\mu(t)},
\end{eqnarray}
where $\mu$ is the probability law of $Y$. When $q = 0$ or in other words $\bbx$ does not have any components, $T(Y, \bbz\mid\bbx) = T(Y, \bbz)$ which measures the dependency of $Y$ on $\bbz$ without any conditioning. 

Based on $n$ samples, the following estimator $T_n$ for $T$ was introduced in \cite{azadkia2021simple}:
\begin{eqnarray}\label{TnContinuous}
T_n(Y, \bbz\mid\bbx) = \frac{\sum_{j=1}^n\min\{R_j, R_{M(j)}\} - \min\{R_j, R_{N(j)}\}}{\sum_{j=1}^n R_j - \min\{R_j, R_{N(j)}\}}.
\end{eqnarray}
Here, $R_j:=\sum_{k=1}^n 1\{Y_k\leq Y_j\}$ is the rank of $Y_j$ among all the other $Y_k$'s, $N(j)$ is the index of the nearest neighbor of $\bbx_j$ among other $\bbx_k$'s, and $M(j)$ is the index of the nearest neighbor of $(\bbx_j, \bbz_j)$ among the rest of $(\bbx_k, \bbz_k)$'s. The nearest neighbors are measured with respect to a metric of choice; in this paper, we use the euclidean distance. Since the estimator~\ref{TnContinuous} only relies on the ranks and the nearest neighbor indices, it is computationally efficient and estimated with computational complexity $O(n\log n)$. For the target variable $Y$, and features $\bbx_1, \cdots, \bbx_p$, FOCI proceeds in the following forward step-wise manner: In each step, FOCI selects one of the remaining variables until none of them adds any ``predictive power". Let $S_t\subseteq[p]$ be the set of selected variables by FOCI up to and including step $t$, initializing with $S_0 = \emptyset$. In step $t+1$, FOCI selects $\argmax_{i\not\in S_t}T_n(Y, \bbx_i\mid\bbx_{S_t})$, if the maximum value is positive. Otherwise FOCI stops and returns $\hat{S} = S_t$. Note that FOCI possesses a natural stopping rule and does not require the user to specify any stopping criterion. We present a summary of FOCI in Algorithm~\oldref{FOCI}.
\begin{algorithm}[h!]
\setstretch{1}
\caption{FOCI~\cite{azadkia2021simple}}\label{FOCI}
\begin{algorithmic}[1]
\Require $n$ samples of $(\bbx_1, \cdots, \bbx_p, Y)$
\Ensure  $\widehat{\MB}(Y)\subseteq [p]$, an estimate of the indices of the Markov blanket of $Y$
\State ${S} = \emptyset$
\While {$\max_{i\not\in{S}} T_n(Y, \bbx_i\mid \bbx_{{S}}) > 0$}
    \State ${S} = {S}\cup \{\argmax_{i\not\in{S}} T_n(Y, \bbx_i\mid \bbx_{{S}})\}$
\EndWhile
\State $\widehat{\MB}(Y) = S$
\end{algorithmic}
\end{algorithm}
In~\cite[Th.6.1]{azadkia2021simple}, the authors show that under mild assumptions, and for large enough sample size $n$, $\widehat{\MB}(Y)$ is with high probability a Markov blanket for $Y$. In general, there is no bound on the size of this Markov blanket, and it may not be the minimal one. In Theorem~\oldref{main} we show that under a so-called $\delta$-Markov Gap assumption~\ref{a3}, the estimated Markov blanket by FOCI is minimal and hence is the Markov boundary. 
\subsection{Stage II of DAG-FOCI: Combining Markov boundaries to identify parents}\label{Parents}

We next describe steps 2 and 3 of Algorithm \oldref{alg:overall}, namely: how DAG-FOCI combines the information of Markov blankets obtained from the first stage to estimate the set(s) of parents. To develop these steps, we assume that the tree-neighborhood assumption for node $Y$ is satisfied (see Definition~\oldref{tnb}); nonetheless, we will show in Theorem~\oldref{FalseDiscovery} that as long as $Y$ is a collider, DAG-FOCI will not output wrong causal conclusions even when the tree-neighborhood assumption does not hold. 

Let $\mathcal{P}(Y) = \{\texttt{pa}_1(Y), \cdots, \texttt{pa}_K(Y)\}$ be the collection of all sets of parents of $Y$ in the Markov equivalence class of $D^\star$ (as defined in Section~\oldref{sec:not}). The logic behind steps 2 and 3 heavily relies on the following properties, satisfied for every parental set $\texttt{pa}_k(Y) \in \mathcal{P}(Y)$:  
\begin{enumerate}[label=(p\arabic*),ref=p\arabic*]
    \item \label{p1} If $|\texttt{pa}_k(Y)|\geq 2$, then for all $i, j\in\texttt{pa}_k(Y)$: $i \in \MB(\bbx_j)$
    
    \item \label{p3} If $\texttt{ch}_k(Y)\neq \emptyset$, then for all $i \in \texttt{pa}_k(Y), j \in \MB(Y){\backslash}{\texttt{pa}_k(Y)}: i \not\in \MB(\bbx_j)$ 
    
    \item \label{p2} For $i, j\in[p]$, $\bbx_i\in\MB(\bbx_j)$, $\bbx_j\in\MB(\bbx_i)$, and $\bbx_i\perp\bbx_j$ iff there exists $k\in[K]$ such that $i, j\in\texttt{pa}_k(Y)$
    
    \item \label{p4} If $|\texttt{pa}_k(Y)| > 1$, we have $K = 1$ and if $K > 1$, we have $|\texttt{pa}_k(Y)| \leq 1$ for all $k\in[K]$
\end{enumerate}
We prove properties \ref{p1}-\ref{p4} in Appendix~\oldref{Appendixp1p6}.

To better understand step 2, we characterize the population analogue of $\hat{\mathcal{S}}_n(Y)$, denoted by $\mathcal{S}^\star(Y)$. The set  $\mathcal{S}^\star(Y)$ is the output of step 2 if the Markov boundary search and pairwise independence tests of DAG-FOCI are performed without making any errors.  Appealing to properties \ref{p1}-\ref{p4}, we conclude that $\mathcal{S}^\star(Y)$ contains all possible sets of parents, i.e. $\mathcal{S}^\star(Y) \supseteq \mathcal{P}(Y)$. Furthermore, $\mathcal{S}^\star(Y)$ contains only nodes that are connected to $Y$ in the skeleton of $D^\star$, i.e. $S\cap \texttt{sp}(Y) = \emptyset$ for every $
S\in \mathcal{S}^\star(Y)$. Thus, in population, step 2 removes all spouses from the Markov boundary of $Y$ without excluding any parental nodes. 

We propose the following procedure to obtain $\mathcal{S}^\star(Y)$ efficiently.  Starting from an empty graph  with the vertex set ${\MB}(Y)$, we form an undirected graph ${G}^\star_Y$ as follows: for every pair of nodes $i, j \in {\MB}(Y)$, edge $i-j$ is added to ${G}^\star_Y$ if $i \in {\MB}(\bbx_j)$ and $j \in {\MB}(\bbx_i)$. Properties~\ref{p1} and~\ref{p3} guarantee that ${G}^\star_Y$ consists of disjoint fully connected subgraphs, i.e. ${G}^\star_Y = G_1\cup\cdots\cup G_R$ such that $G_i$'s are complete graphs and $V(G_i)\cap V(G_j) = \emptyset$. Among these connected components, we let ${\mathcal{S}}^\star(Y)$ be those that contain pairwise independent nodes. Algorithm~\oldref{alg:stage2} presents the finite-sample analogue to obtain the estimate $\hat{\mathcal{S}}_n(Y)$.

\begin{algorithm}[h]
\setstretch{1}
\caption{Stage II of DAG-FOCI (implementation of steps 2 in Algorithm \oldref{alg:overall})}
\begin{algorithmic}[1]
\Require $\widehat{\MB}(Y)$ and $\widehat{\MB}(\bbx_i)$ for all $i\in\widehat{\MB}(Y)$
\Ensure $\hat{\mathcal{S}}_n(Y)$
\State $\hat{\mathcal{S}}_n(Y) = \{\emptyset\}$
\State Create ${G}^\star_Y = (\tilde{V}=\widehat{\MB}(Y), \tilde{E})$ where $(i, j)\in\tilde{E}$ iff $i\in\widehat{\MB}(\bbx_j)$ and $j\in\widehat{\MB}(\bbx_i)$.

\State Decompose ${G}^\star_Y$ to its disjoint connected components $G_1\cup\cdots\cup G_R$. 

\State For all $r\in[R]$ add the set of vertices of $G_r$ to $\hat{\mathcal{S}}_n$ iff for any pair $i, j$ in $G_r$, $\bbx_i\perp\bbx_j$.
\end{algorithmic}
\label{alg:stage2}
\end{algorithm}

Thus far, we have presented an efficient approach to combine local Markov structure with pairwise independence tests to arrive at a superset $\mathcal{S}^\star(Y)$ of the true parental set $\mathcal{P}(Y)$. We now describe step 3, which removes non-parental elements of $\mathcal{S}^\star(Y)$ when possible. To that end, we highlight the following relation $\mathcal{S}^\star(Y) =  \mathcal{P}(Y) \cup \mathcal{E}$ where $\mathcal{E}$ satisfies:
\begin{enumerate}[label=(e\arabic*),ref=e\arabic*]
    \item\label{e1} $\cup_{k=1}^K\texttt{pa}_k(Y)\cap S = \emptyset$ for all $S\in\mathcal{E}$
    \item\label{e2} $|S| \leq 1$ for all $S\in\mathcal{E}$
\end{enumerate}
Properties \ref{e1} and \ref{e2} are implied by \ref{p1}-\ref{p4} and motivate step 3 of Algorithm~\oldref{alg:overall}. Specifically, according to property \ref{e1}, the set $\mathcal{E}$ consists of non-parental elements and should ideally be removed from $\mathcal{S}^\star(Y)$. Property \ref{e2} guides how we remove $\mathcal{E}$ (when possible) in step 3.

Specifically, consider all the possible scenarios for the set $\mathcal{S}^\star(Y)$: $(i)$ there exists a single set $S \in \mathcal{S}^\star(Y)$ with $|S|>1$, $(ii)$ there exists multiple sets $S \in \mathcal{S}^\star(Y)$ with $|S| > 1$, and $(iii)$ $|S| \leq 1$ for all $S \in \mathcal{S}^\star(Y)$. Appealing to properties \ref{e1}-\ref{e2} and \ref{p4} , scenario $(i)$ implies that $\mathcal{P}(Y) = \{S\} = \{\texttt{pa}(Y)\}$ and scenario $(ii)$ leads to a contradiction. These observations explain the following logic in step 3 of Algorithm~\oldref{alg:overall}: if $\hat{\mathcal{S}}_n(Y)$ contains exactly one set $S$ such that $|S|>1$, we let $\hat{\mathcal{P}}_n(Y)= \{S\}$; if  there are multiple sets $S$ with $|S|>1$, DAG-FOCI returns an error message. Finally, in scenario $(iii)$, we cannot distinguish  $\mathcal{E}$ from $\mathcal{P}(Y)$, as property \ref{p4} states that $\mathcal{P}(Y)$ can consist of all singletons. This observation explains the following logic in step 3 of Algorithm~\oldref{alg:overall}: if $|\mathcal{S}| \leq 1$ for all $S \in\hat{\mathcal{S}}_n(Y)$, we let $\hat{\mathcal{P}}_n(Y) = \hat{\mathcal{S}}_n(Y)$.

\begin{figure}[ht]
\centering
  \begin{tikzpicture}
        [->, node/.style={}]
        \node[node]    (X1)                     {$\bbx_1$};
        \node[node]    (X4)    [right= of X1]   {$\bbx_4$};
        \node[node]    (X5)    [left= of X1]    {$\bbx_5$};
        \node[node]    (Y)     [below= of X1]   {$Y$};
        \node[node]    (X2)    [right= of Y]    {$\bbx_2$};
        \node[node]    (X3)    [left= of  Y]    {$\bbx_3$};
        
        \path (X1) edge (Y);
        \path (Y) edge (X2);
        \path (Y) edge (X3);
        \path (X4) edge (X1);
        \path (X5) edge (X1);
    \end{tikzpicture}
    \caption{DAG $D$ used to illustrate the case, where the output of DAG-FOCI is a superset of the set of all possible parental sets of $Y$}
    \label{fig:nonmultipleparent}
\end{figure}
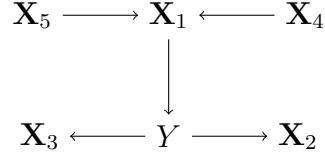

We remark that scenario $(iii)$ highlights a fundamental challenge with applying a local algorithm, like DAG-FOCI, to obtain the parental set(s) of $Y$. Specifically, let $D^\star$ be the DAG in Figure~\oldref{fig:nonmultipleparent}. Even if all of the steps (e.g. the Markov boundary search and pairwise independence tests) of DAG-FOCI are performed without making an error, DAG-FOCI will output  $\hat{\mathcal{P}}_n(Y) = \{\{1\}, \{2\}, \{3\}, \emptyset\}$ as the collection of parental sets, whereas the true parental set is $\{3\}$. This simple example demonstrates that when $Y$ has only a single parent, DAG-FOCI may produce a superset of the true parental set. 
To overcome this problem, one can modify the result of DAG-FOCI by running DAG-FOCI over neighborhoods of $Y$ and combine the outputs. For example in $D$ in Figure~\oldref{fig:nonmultipleparent}, DAG-FOCI with high probability gives $\hat{\mathcal{P}}_n(\bbx_1) = \{\{4, 5\}\}$ which rules out the possibility of $\bbx_1$ being a child of $Y$. Hence we can modify the output of DAG-FOCI for $Y$ and remove the extra sets $\{2\}$ and $\{3\}$ from its output. In general, this process may require investigating a larger neighborhood of $Y$ and cannot be limited to only members of $\MB(Y)$.

\section{Theoretical Guarantees for DAG-FOCI}
\label{sec:Theory}

Let $\bbx = (\bbx_1, \cdots, \bbx_p, Y = \bbx_{p+1})$ be as in the previous section, where $\bbx_i$'s can be continuous or discrete or a mix of both. For any nonempty set $S\subseteq[p]$, consider the function $Q$ defined as in (6.1) in~\cite{azadkia2021simple}, where $Q(\emptyset) = 0$ and 
\begin{eqnarray}
Q(S) := Q(Y, \bbx_S) = \int\var(\pp(Y\geq t\mid \bbx_S))d\mu(t).
\end{eqnarray}
Note that $Q(Y, \bbx_S)$ is the numerator of $T(Y, \bbx_S)$ in~\ref{TOriginal}.

Let $\delta$ be the largest number such that for any insufficient subset $S$ for $Y$ ($S$ not sufficient), there is some $i\not\in S$ such that $Q(S\cup\{i\})\geq Q(S) + \delta$. In other words adding $i$ to $S$ increases the ``predictive power" by at least $\delta$. The definition of $\delta$ ensures that there is at least one sufficient subset of size at most $1/\delta$.

Our theoretical guarantees require the following three assumptions on the joint distribution of $\bbx$:
\begin{enumerate}[label=(A.\arabic*),ref=A.\arabic*]
    \item \label{a1} There are non-negative real numbers $\beta$ and $C$ such that for any set $S\subseteq\{1, \cdots, p\}$ of size $\leq 1/\delta + 2$, any $x, x^\prime\in\rr^S$ and any $t\in\rr$, 
    \begin{eqnarray*}
    |\pp(Y\geq t\mid\bbx_S = x) - \pp(Y\geq t\mid\bbx_S = x^\prime)| \leq C(1 + \|x\|^\beta + \|x^\prime\|^\beta)\|x - x^\prime\|.
    \end{eqnarray*}
    \item \label{a2} There are positive numbers $C_1$ and $C_2$ such that for any $S$ of size $\leq 1/\delta + 2$ and any $t > 0$, $\pp(\|\bbx_S\|\geq t)\leq C_1e^{-C_2 t}$.
    \item \label{a3} ($\delta$-Markov Gap) For any subset $S\subseteq [p]$ such that $S\subsetneq\MB(Y)$, there exists $i$ with $i\in\MB(Y)\setminus S$, such that for any $j$ with $j\not\in\MB(Y)$
    \begin{eqnarray*}
    Q(S\cup\{i\}) - Q(S\cup\{j\})\geq\delta/4.
    \end{eqnarray*}
\end{enumerate}
The value $\delta > 0$ is the same everywhere in this paper for simplicity of notation. Furthermore, $\|\cdot\|$ denotes the euclidean norm. Assumptions~\ref{a1} and~\ref{a2} are identical to the ones in \cite{azadkia2021simple} for identifying a Markov blanket. Assumption \ref{a3}, called $\delta$-Markov Gap, ensures that the estimated Markov blanket is minimal. More specifically,~\ref{a3} requires that the members of the Markov boundary have higher predictive power compared to the rest of the variables.

In~\cite[Th.6.1]{azadkia2021simple}, under Assumptions~\ref{a1}, and~\ref{a2}, the authors show that the output  $\widehat{\MB}(Y)$ of FOCI is with high probability a Markov blanket for $Y$. With the addition of Assumption~\ref{a3}, we show that FOCI outputs the Markov boundary (i.e. smallest Markov blanket):
\begin{thm}\label{main}
Suppose Assumptions~\ref{a1}-\ref{a3} hold for some (and the same) $\delta > 0$. Let $\widehat{\MB}(Y)$ be the subset selected by FOCI (Algorithm~\oldref{FOCI}) with a sample of size $n$. There are positive real numbers $L_1$, $L_2$, and $L_3$ depending only on $C$, $\beta$, $C_1$, $C_2$, and $\delta$ such that $\pp(\widehat{\MB}(Y) = \MB(Y))\geq 1 - L_1 p^{L_2}e^{-L_3 n}$.
\end{thm}
In the following lemma, we show how the $\delta$-Markov Gap is satisfied for some $\delta^\prime>0$ when $Y$ has the tree-neighborhood property (see Definition \oldref{tnb}):
\begin{lmm}\label{prop1}
Assume there exists $\delta > 0$ such that for any insufficient set $S$ there exists $i$ such that $Q(S\cup\{i\})\geq Q(S) + \delta$. If $Y$ has the tree-neighborhood property then there exists $0 < \delta^\prime\leq\delta$ such that $Y$ satisfies the $\delta$-Markov Gap property in~\ref{a3} with $\delta^\prime$.
\end{lmm}
Together with Theorem~\oldref{main}, Lemma~\oldref{prop1} implies that for any node with the tree-neighborhood property, FOCI outputs their Markov boundary with high probability.

So far, we have analyzed the first stage of DAG-FOCI. In order to extend our analysis to the entire DAG-FOCI procedure, we must assess the accuracy of the independence tests in the second stage of DAG-FOCI. To that end, we introduce some notations. Given $n$ i.i.d observations $\{\bbv_i, \bbw_i\}_{i=1}^n$ of real-valued random variables $\bbw$ and $\bbv$, a statistical test $\phi_\alpha:(\rr\times\rr)^n\rightarrow\{0, 1\}$ distinguishes between the null hypothesis $\mathcal{H}_0:\bbw\perp\bbv$ and the alternate hypothesis $\mathcal{H}_1:\bbw\not\perp\bbv$. Here, $\alpha \in [0,1]$ is the size of the test and is an upper bound on the Type-I error, i.e. the probability that the test $\phi_\alpha$ falsely rejects $\mathcal{H}_0$. Associated with the test $\phi_\alpha$ is also the Type-II error $\beta_{\bbw, \bbv}^{\phi_\alpha}(n) \in [0,1]$, which is the probability that $\phi_\alpha$ fails to reject the null when $\bbw\not\perp\bbv$. In our analysis, we fix $\alpha$ and the test $\phi_\alpha$, and define the quantity
\begin{eqnarray}
\beta(n) := \max_{i, j\in\MB(Y)\text{ s.t. }\bbx_i\not\perp\bbx_j} \beta_{\bbx_i, \bbx_j}^{\phi_{\alpha}}(n).
\label{betan}
\end{eqnarray}
Here, $\beta(n)$ is the maximum Type-II error among all tests between pairs of dependent random variables in the Markov boundary of $Y$. 

Now that we have characterized the performance of pairwise independence independence tests in the second stage of DAG-FOCI, we provide finite-sample guarantees for the entire procedure of DAG-FOCI.
\begin{thm}\label{bonferroni}
Suppose Assumptions \ref{a1}-\ref{a3} hold for some $\delta>0$ and $Y$ satisfies the tree-neighborhood property. Then:
\begin{eqnarray*}
\pp(\mathcal{P}(Y)\not\subseteq\hat{\mathcal{P}}_n(Y)) &\leq & (|\MB(Y)|+1)(L_1 p^{L_2}e^{-L_3 n}) + (|\MB(Y)| - 1)\beta(n) + \\ && \alpha\max_k|\texttt{pa}_k(Y)|(\max_k|\texttt{pa}_k(Y)| - 1)/2.
\end{eqnarray*}
Constant $\alpha$ is the size of the independence tests in Algorithm~\oldref{alg:overall} and $\beta(n)$ is defined in \ref{betan}. The positive numbers $L_1$, $L_2$, and $L_3$ are depending only on $C$, $\beta$, $C_1$, $C_2$, and $\delta$. 
\end{thm}
Theorem~\oldref{bonferroni} guarantees that for large enough $n$ and with high probability, DAG-FOCI's output contains the true parental sets. We remark here that while $\hat{\mathcal{P}}_n(Y)$ may be a superset of $\mathcal{P}(Y)$, it can still provide meaningful causal conclusions. Specifically, it is straightforward to show that with high probability, $\hat{\mathcal{P}}_n(Y)$ contains only nodes that are directly connected to $Y$ in the skeleton of $D^\star$. Nonetheless, the true parental sets $\mathcal{P}(Y)$ can be recovered under the tree-neighborhood assumption and if $Y$ is a collider (i.e. has multiple parents):
\begin{cor}\label{corbon}
With the same setting and constants as in Theorem~\oldref{bonferroni} and in addition when $\mathcal{P}(Y) = \{\texttt{pa}(Y)\}$ and $|\texttt{pa}(Y)|>1$ we have
\begin{eqnarray*}
\pp(\hat{\mathcal{P}}_n(Y)\neq\mathcal{P}(Y)) &\leq & (|\MB(Y)|+1)(L_1 p^{L_2}e^{-L_3 n}) + (|\MB(Y)| - 1)\beta(n) + \\ && \alpha\max_k|\texttt{pa}_k(Y)|(\max_k|\texttt{pa}_k(Y)| - 1)/2.
\end{eqnarray*}
\end{cor}
Note that the probability bound in Theorem~\oldref{bonferroni} depends on some values that are usually unknown, namely $|\MB(Y)|$, and $|\texttt{pa}_k(Y)|$. They can be bounded as follows. Denote by $d$ the maximal degree of the skeleton of $D^\star$. Then we have: $(i)$ $|\MB(Y)|\leq d(d-1)$ as $Y$ has at most $d$ directly connected neighbors and each of those have at most $d-1$ nodes connected to them besides $Y$; $(ii)$ $|\texttt{pa}_k(Y)|\leq d$; Putting these together shows that the output of DAG-FOCI contains the sets of parents of $Y$ with probability at least $1 - d^2\max (L_1 p^{L_2}e^{-L_3 n}, \alpha, \beta(n))$. Thus we guarantee that DAG-FOCI captures all the sets of parents correctly with high probability, even for high-dimensional cases where $p$ can grow polynomially fast with the sample size but the graph is sparse with maximal degree $d$ being constant or growing sufficiently slowly with $n$.

\subsection{Beyond polytrees}
The tree-neighborhood assumption may not be satisfied in practice. In the following theorem, we characterize the output of DAG-FOCI without the tree-neighborhood assumption, as long as $Y$ is an unshielded collider: 
\begin{thm}\label{FalseDiscovery}
Suppose~\ref{a1}-\ref{a3} hold for $Y$ for some $\delta > 0$ and~\ref{a1}-\ref{a2} hold for all members of $\MB(Y)$. Let $\hat{\mathcal{P}}_n(Y)$ be the output of DAG-FOCI for $Y$. Assume that $\mathcal{P}(Y) = \{\texttt{pa}(Y)\}$ and $|\texttt{pa}(Y)|>1$. Let $F$ be the event that $\exists S\in\hat{\mathcal{P}}_n(Y)$ such that $S\setminus \texttt{pa}(Y)\neq\emptyset$. Then
\begin{eqnarray*}
\pp(F) \leq \frac{1}{2}|\MB(Y)|^2(\beta(n) + \alpha) + (|\MB(Y)| + 1)L_1p^{L_2}e^{L_3n}.
\end{eqnarray*}
The positive real numbers $L_1$, $L_2$, and $L_3$ are depending only on $C$, $\beta$, $C_1$, $C_2$, and $\delta$. Constant $\beta(n)$ is defined in \ref{betan} and $\alpha > 0$ is the size of the independence tests. 
\end{thm}
Theorem~\oldref{FalseDiscovery} emphasizes the conservative nature of DAG-FOCI: when the data does not meet our assumptions, then with high probability, DAG-FOCI does not yield any false discoveries.


\section{Extensions of DAG-FOCI to interventional data}\label{sec:interventional}
In Section~\oldref{sec:DAG-FOCI}, we developed DAG-FOCI in a general observational setting where the parental sets of a target variable can only be identified up to the Markov equivalence class. In the general case, however, identifiability can only be improved by performing interventions (experiments). Examples of such interventions are abundant in the empirical sciences, including in the dataset we analyze in Section~\oldref{Real Data} on protein expressions. Given the potential benefits of interventions and their applicability in real data, we develop an extension of DAG-FOCI to exploit interventional data. 

We consider a setting where we have access to an observational environment and at least a single interventional environment with a do-intervention on $Y$ (interventions on other variables are also allowed). Once again we assume that $Y$ satisfies the tree-neighborhood assumption in the population DAG $D^\star$. We first supply data from the observational environment to DAG-FOCI (in Algorithm~\oldref{alg:overall}) to obtain the set of all possible parents $\hat{\mathcal{P}}_n(Y)$; from Theorem~\oldref{bonferroni}, we conclude that $\hat{\mathcal{P}}_n(Y)$ is a superset of the true parental sets $\mathcal{P}(Y)$ with high probability. We then apply FOCI (in Algorithm~\oldref{FOCI}) to the interventional environment to obtain a Markov boundary $\widehat{\MB}_I(Y)$ of $Y$; from Theorem~\oldref{main}, we conclude that $\widehat{\MB}_I(Y)$ coincides with the true Markov boundary of $Y$ in the modified DAG after interventions. Appealing to the fact that a do-intervention on $Y$ cuts off its dependencies with its parents, as well to the tree-neighborhood property, we conclude that $\widehat{\MB}_I(Y)$ cannot contain any parents from the population DAG $D^\star$. Thus, for any $S \in \hat{\mathcal{P}}_n(Y)$ where $S \cap \widehat{\MB}_I(Y) \neq \emptyset$, the set of variables $S \cap \widehat{\MB}_I(Y)$ are child nodes of $Y$ in the population DAG $D^\star$. Furthermore, again due to the tree-neighborhood assumption, such a set $S$ cannot contain a parent of $Y$ in $D^\star$. Thus, following the recipe outlined above, which is also summarized in Algorithm~\oldref{alg:interven}, we combine observational and interventional data to arrive at a modified parental set $\tilde{\mathcal{P}}_n(Y)$ and child nodes $\tilde{C}_n(Y)$ for the target variable $Y$.
\begin{algorithm}[h]
\setstretch{1}
\caption{DAG-FOCI with interventional data}
\label{alg:interven}
\begin{algorithmic}[1]
\Require samples of $(\bbx_1, \cdots, \bbx_p, Y = \bbx_{p+1})$ in an observational environment and samples of $(\bbx_1, \cdots, \bbx_p, Y = \bbx_{p+1})$ in an interventional environment with a do-intervention on $Y$
\Ensure  $\tilde{\mathcal{P}}_n(Y)$: the modified set of all possible parental sets after combining observational and interventional data; $\tilde{C}_n(Y)$: a set containing (some) children of $Y$
\State {\bf DAG-FOCI on observational data}: supply observational data to DAG-FOCI (Algorithm~\oldref{alg:overall}) to obtain all possible parental sets $\hat{\mathcal{P}}_n(Y)$
\State {\bf Markov boundary of $Y$ in interventional environment}: supply interventional data to FOCI (Algorithm~\oldref{FOCI}) to obtain $\widehat{\MB}_I(Y)$, the Markov boundary of $Y$ in the interventional environment
\State {\bf Using interventional information to identify parent and child nodes}: initialize $\tilde{\mathcal{P}}_n(Y) = \emptyset$ and $\tilde{C}_n(Y) = \emptyset$ and perform for all $S \in \hat{\mathcal{P}}_n(Y)$:
\begin{itemize} \item if $S \cap \widehat{\MB}_I(Y) = \emptyset$, add $S$ to $\tilde{\mathcal{P}}_n(Y)$. Otherwise, add $S  \cap \widehat{\MB}_I(Y)$ to $\tilde{C}_n(Y)$. 
\end{itemize}
\end{algorithmic}
\end{algorithm}
As a straightforward result of Theorems~\oldref{main} and \oldref{bonferroni}, we next provide theoretical guarantees that the output of Algorithm~\oldref{alg:interven} uniquely identifies (some of) the children and parents of node $Y$ in the population DAG $D^\star$. Letting $\texttt{pa}(Y)$ and $\texttt{ch}(Y)$ be the parents and children of $Y$ in the population DAG $D^\star$, we have the following theorem:

\begin{cor}\label{thm:interv}
Suppose that Assumptions~\ref{a1}-\ref{a3} hold for some $\delta>0$ for both the observational and interventional environments. Then outputs $\tilde{C}_n(Y)$ and $\tilde{\mathcal{P}}_n(Y)$ of Algorithm~\oldref{alg:interven} satisfy
\begin{eqnarray*}
\mathbb{P}(\tilde{C}_n(Y) \not\subseteq \texttt{ch}(Y)) &\leq & (|\MB(Y)|+1)(L_1 p^{L_2}e^{-L_3 n_\text{obs}}) +\\
&& |\MB(Y)|^2\beta(n_\text{obs}) + L_1p^{L_2}e^{-L_3 n_I}
\end{eqnarray*}
Furthermore, assuming that there are no do-interventions on any child node of $Y$, that $Y$ satisfies the tree-neighborhood property in Definition~\oldref{tnb}, we have:
\begin{eqnarray*}
\mathbb{P}(\tilde{C}_n(Y) \not\subseteq \texttt{ch}(Y) \cup \tilde{\mathcal{P}}_n(Y) \neq \{\texttt{pa}(Y)\}) &\leq & (|\MB(Y)|+1)(L_1 p^{L_2}e^{-L_3 n_\text{obs}}) +\\
&& (|\MB(Y)| - 1)\beta(n_\text{obs}) +\\ 
&& \alpha\max|\texttt{pa}_k(Y)|(\max|\texttt{pa}_k(Y)| - 1)/2 + \\
&& L_1p^{L_2}e^{-L_3 n_I},
\end{eqnarray*}
where $n_\text{obs}$ and $n_I$ are the number of samples in the observational and interventional environment, respectively and all the constants are as in Theorem~\oldref{bonferroni}. 
\label{corr:interv}
\end{cor}
Corollary~\oldref{thm:interv} ensures, with high probability, that any child $\tilde{C}_n(Y)$ identified by Algorithm~\oldref{alg:interven} is a child of $Y$ in the population DAG $D^\star$. Note that this result does not require $Y$ to satisfy the tree-neighborhood assumption, and also does not rule out the possibility of having interventions on the children. if there are no interventions on any child nodes of $Y$ and under the tree-neighborhood assumption on $Y$, Corollary~\oldref{thm:interv} ensures that with high probability, $\tilde{\mathcal{P}}_n(Y)$ is the parental set of $Y$ in $D^\star$. 

\section{Empirical results}\label{sec:Experiments}
In this section, we study the performance of the DAG-FOCI algorithm on both synthetic and real data examples. We compare it with Hill climbing~\cite{gamez2011learning}, two version of the PC-algorithm, namely for the Gaussian case with partial correlations and for the non-parametric case with HSIC~\cite{fukumizu2007kernel} for conditional independence testing, CAM~\cite{buhlmann2014cam}, and GES~\cite{chickering2002optimal}. 

In the second stage of DAG-FOCI we perform all the independence tests using permutation test and measure of dependence CODEC~\cite{azadkia2021simple}. We generate 100 independent permutations for each test over the set of integers $[n]$, where $n$ is the sample size, and we set our threshold to $\alpha = 0.05$.
\subsection{Simulated Data}\label{Simulated Data}

All the results are based on 100 independent simulation runs. 
\begin{ex}\label{example05}
We consider the DAG structure in Figure~\oldref{fig:example05} with the corresponding relationship between the nodes described in the following set of Equations~\ref{eq_ex_05}. 
\begin{equation}\label{eq_ex_05}
\begin{split}
&\bbx_1, \bbx_2, \bbx_3, \bbx_4, \bbx_8, \bbx_{10},\bbx_{12}, \varepsilon_i \overset{\text{i.i.d}}{\sim}\mathcal{N}(0, 1), \\
&\bbx_5 = \bbx_1 - \arctan(\bbx_2) + \varepsilon_5, \\
&\bbx_6 = \bbx_2 + \bbx_4 + \bbx_3^2 + \varepsilon_6, \\
&\bbx_7 = \sin(\bbx_3) + \varepsilon_7, \\
&\bbx_9 = \sin(\bbx_6 + \varepsilon_9) + |\bbx_{10}|, \\
\end{split}
\qquad
\begin{split}
&\bbx_{11} = \bbx_6(\bbx_{12} - \bbx_{8}) + \varepsilon_{11}, \\
&\bbx_{13} = \arctan(\bbx_9^2 + \varepsilon_{13}),\\
&\bbx_{14} = \sin(\bbx_{11}) + \varepsilon_{14}, \\
&\bbx_{15} = \sqrt{|\bbx_{12}|} + \varepsilon_{15}, \\
&\bbx_{16} = \sin(\bbx_{12}) + \varepsilon_{16}. 
\end{split}
\end{equation}

\begin{figure}[h!]
\centering
\begin{tikzpicture}
[->, node/.style={}]
\node[node]    (X1)                     {$\bbx_1$};
\node[node]    (X2)    [right= of X1]   {$\bbx_2$};
\node[node]    (X3)    [right= of X2]   {$\bbx_3$};
\node[node]    (X4)    [right= of X3]   {$\bbx_4$};
\node[node]    (X5)    [below= of X1]   {$\bbx_5$};
\node[node]    (X6)     [below= of X2]  {$\bbx_6$};
\node[node]    (X7)    [below= of X3]   {$\bbx_7$};
\node[node]    (X8)    [below= of X4]   {$\bbx_8$};
\node[node]    (X9)    [below= of X5]   {$\bbx_9$};
\node[node]    (X10)    [below= of X6]   {$\bbx_{10}$};
\node[node]    (X11)    [below= of X7]   {$\bbx_{11}$};
\node[node]    (X12)    [below= of X8]   {$\bbx_{12}$};
\node[node]    (X13)    [below= of X9]   {$\bbx_{13}$};
\node[node]    (X14)    [below= of X10]   {$\bbx_{14}$};
\node[node]    (X15)    [below= of X11]   {$\bbx_{15}$};
\node[node]    (X16)    [below= of X12]   {$\bbx_{16}$};

\path (X1) edge (X5);
\path (X2) edge (X5);
\path (X2) edge (X6);
\path (X3) edge (X6);
\path (X3) edge (X7);
\path (X4) edge (X6);
\path (X6) edge (X9);
\path (X6) edge (X11);
\path (X8) edge (X11);
\path (X9) edge (X13);
\path (X10) edge (X9);
\path (X11) edge (X14);
\path (X12) edge (X11);
\path (X12) edge (X15);
\path (X12) edge (X16);
\end{tikzpicture}
\caption{DAG in Setting~\ref{example05}. The goal is to infer the parents of $\bbx_6$ and $\bbx_{11}$.} 
\label{fig:example05}
\end{figure}
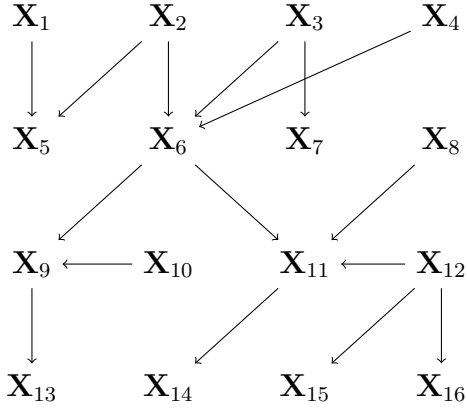

We study the performance of DAG-FOCI in estimating the set of parents of $\bbx_6$ and $\bbx_{11}$, which  are both identifiable,  over 100 independent simulation runs for sample sizes $n = 2000, 4000, 6000, 8000, 10000$.
 
We evaluate our estimations using the Jaccard index. The Jaccard (similarity) index between two sets $A$ and $B$ is defined as $|A\cap B|/|A\cup B|$ and therefore, it takes values in $[0, 1]$.
\begin{figure}[h!]
  \centering
  \includegraphics[width=1\textwidth]{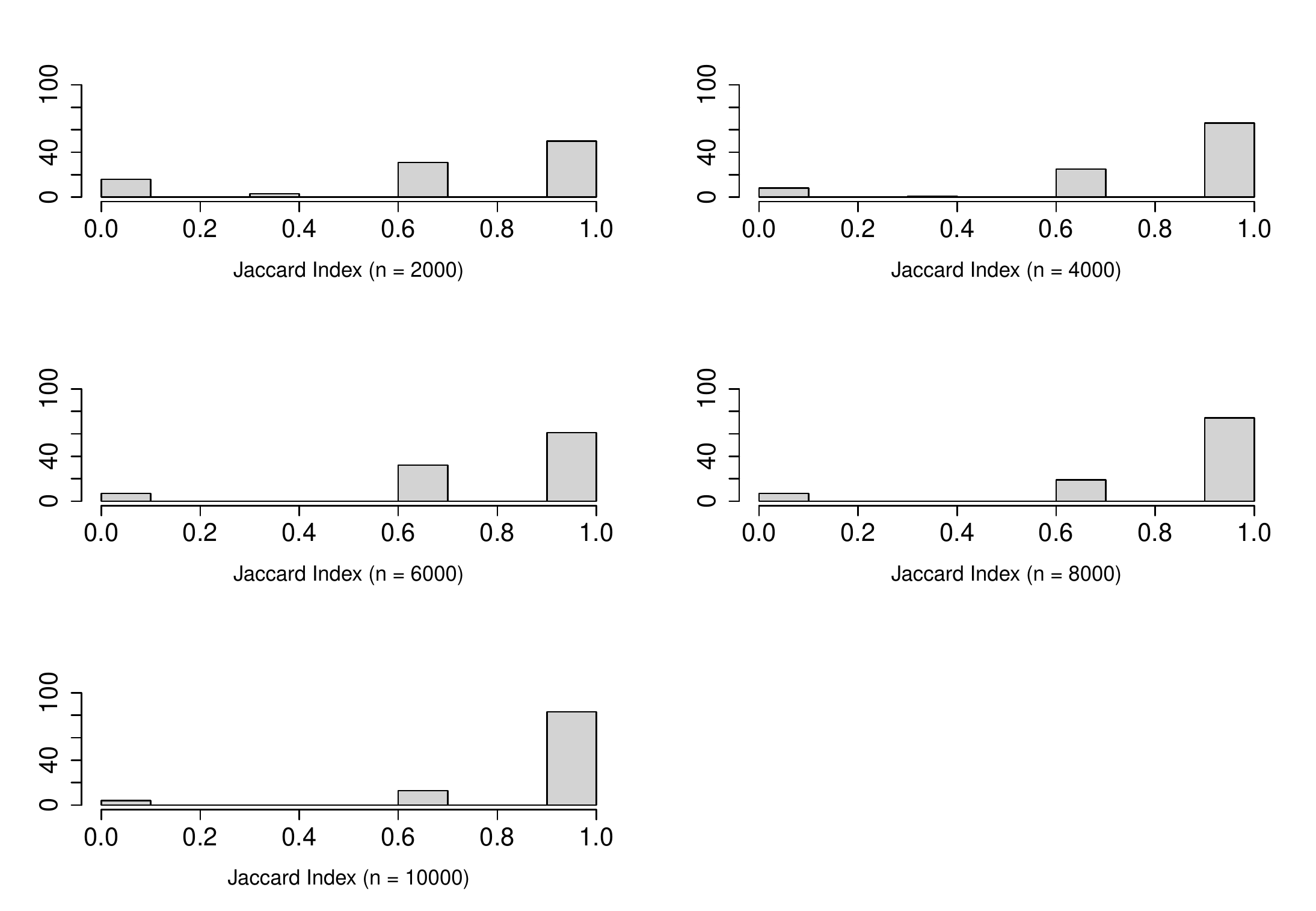}
    \caption{Histogram of the Jaccard index between the estimated set of parents and true parent set of node $\bbx_6$ in Setting~\ref{example05} for sample sizes $n = 2000, 4000, 6000, 8000, 10000$.}
    \label{JaccardX6ex05}
\end{figure}

\begin{figure}[h!]
  \centering
  \includegraphics[width=1\textwidth]{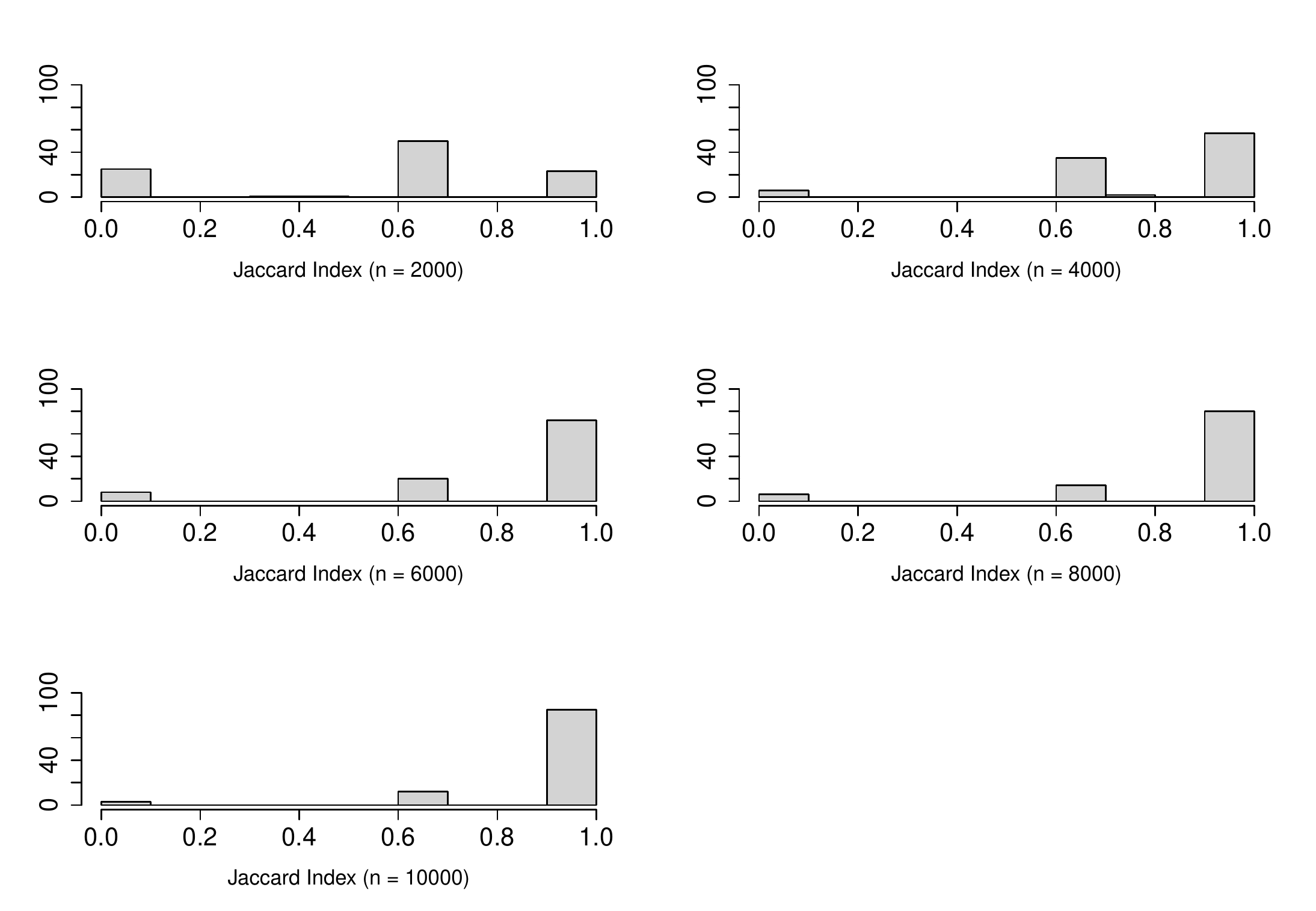}
    \caption{Histogram of the Jaccard index between the estimated set of parents and true parent set of node $\bbx_{11}$ in Setting~\ref{example05} for sample sizes $n = 2000, 4000, 6000, 8000, 10000$.}
    \label{JaccardX11ex05}
\end{figure}
Figures~\oldref{JaccardX6ex05} and~\oldref{JaccardX11ex05} show the histograms of the Jaccard index between the estimated set of parents and the true set of parents for nodes $\bbx_6$ and $\bbx_{11}$ respectively. Note that $\texttt{pa}(\bbx_6) = \{2, 3, 4\}$ and $\texttt{pa}(\bbx_{11}) = \{6, 8, 12\}$. If $\hat{\mathcal{P}}_n(\bbx_6)$ and $\hat{\mathcal{P}}_n(\bbx_{11})$ contains more than one set, we set the Jaccard index of that run equal to zero. As the sample size $n$ grows, the Jaccard index concentrates more on $1$, which is expected from Theorem~\oldref{bonferroni}. Table~\oldref{TableSummaryEx01} summarizes the results of these simulations for different sample sizes.

\begin{table}[h!]
\begin{center}
\centering
\def\arraystretch{1}
\begin{tabular}{lrrrrr}
  \hline
 $\bbx_6$ & Exact recovery & Non-unique & False & Missing & Jaccard Index\\ 
  \hline
    $n = 2000$ & 50 & 13 & 0.01 & 0.85 & 0.72\\
    $n = 4000$ & 66 & 3 & 0.03 & 0.51 & 0.83\\
    $n = 6000$ & 61 & 2 & 0.03 & 0.53 & 0.82\\
    $n = 8000$ & 74 & 0 & 0.05 & 0.4 & 0.87\\
    $n = 10000$ & 83 & 0 & 0.02 & 0.25 & 0.92\\
   \hline
   $\bbx_{11}$ & Exact recovery & Non-unique & False & Missing & Jaccard Index\\ 
  \hline
    $n = 2000$ & 23 & 18 & 0.08 & 1.28 & 0.57\\
    $n = 4000$ & 57 & 2 & 0.06 & 0.53 & 0.82\\
    $n = 6000$ & 72 & 3 & 0.05 & 0.44 & 0.85\\
    $n = 8000$ & 80 & 2 & 0.04 & 0.32 & 0.90\\
    $n = 10000$ & 85 & 0 & 0.03 & 0.21 & 0.93\\
   \hline
\end{tabular}
\end{center}
\caption{Summary of performance of DAG-FOCI for $\bbx_6$ and $\bbx_{11}$ over 100 simulation runs. Column ``Exact recovery" shows the number of runs in which DAG-FOCI has returned exactly the true set of parents. Column ``Non-unique" shows the number of runs in which DAG-FOCI returns multiple sets. Column ``False" shows the the average number of falsely discovered parents. Column ``Missing" shows the average number of missing parents. Column ``Jaccard Index" shows the average Jaccard index of all 
runs.}
\label{TableSummaryEx01}
\end{table}

\paragraph{Comparison with other methods.}
For sample size $n = 10^4$ and over the 100 independent simulation runs,
we compare the performance of DAG-FOCI in identifying the set of parents of $\bbx_6$ and $\bbx_{11}$ with other methods. We use the implementation of Hill climbing from the R package \texttt{bnlearn}~\cite{bnlearn}, and PC with partial correlation and GES from R package \texttt{pcalg}~\cite{kalisch2010pcalg}, and CAM from R package \texttt{CAM}~\cite{CAM}, and PC with HSIC from R package \texttt{kpcalg}~\cite{kpcalg}. In Table~\oldref{TableComparisonEx01} we summarize the performance of these methods. 
\begin{table}[h!]
\begin{center}
\centering
\def\arraystretch{1}
\begin{tabular}{lrrrrr}
  \hline
 $\bbx_6, (n = 10^4)$ & Exact recovery & Non-unique & False & Missing & Jaccard Index\\ 
  \hline
    DAG-FOCI & 83 & 0 & 0.02 & 0.25 & 0.92\\
    PC (partial correlation) & 3 & 0 & 0.34 & 0.87 & 0.65\\
    PC (HSIC) & \textbf{100} & 0 & \textbf{0} & \textbf{0} & \textbf{1}\\
    Hill-climbing & 4 & 0 & 0.85 & 0.93 & 0.55\\
    CAM & 2 & 0 & \textbf{0} & 2.07 & 0.61\\
    GES & 4 & 0 & 0.88 & 0.94 & 0.54\\
   \hline
   $\bbx_{11}, (n = 10^4)$ & Exact recovery & Non-unique & False & Missing & Jaccard Index\\ 
  \hline
    DAG-FOCI & \textbf{85} & 0 & \textbf{0.03} & 0.21 & \textbf{0.93}\\
    PC (partial correlation) & 1 & 0 & 1.07 & 0.91 & 0.52\\
    PC (HSIC) & 0 & 0 & 1.84 & \textbf{0.09} & 0.60\\
    Hill-climbing & 5 & 0 & 0.49 & 0.94 & 0.60\\
    CAM & 1 & 0 & 0.48 & 2.24 & 0.22\\
    GES & 3 & 0 & 0.59 & 0.9 & 0.60\\
   \hline
\end{tabular}
\end{center}
\caption{Comparison of the result different algorithm for estimating the set of parents of $\bbx_6$ and $\bbx_{11}$ over 100 simulation runs. Column ``Exact recovery" shows the number of runs in which each method has returned exactly the true set of parents. Column ``Non-unique" shows the number of runs in which each method returns multiple sets. Column ``False" shows the the average number of falsely discovered parents. Column ``Missing" shows the average number of missing parents. Column ``Jaccard Index" shows the average Jaccard index of all runs.}
\label{TableComparisonEx01}
\end{table}
As Table~\oldref{TableComparisonEx01} shows, for $\bbx_{11}$ DAG-FOCI has the best results in estimating the exact set of parents. For $\bbx_6$, PC algorithm using HSIC works perfectly in all the runs. Note that the signal in this case is additive and in for $\bbx_{11}$ there are multiplicative interaction terms. This shows that the choice of the kernel plays an important role for the PC algorithm with a Kernel test of conditional independence even with a sample size as large as $10^4$. Although the performance of PC algorithm for $\bbx_6$ is impressive but the same technique fails for $\bbx_{11}$. The low average of false-positive discoveries for DAG-FOCI in Tables~\oldref{TableSummaryEx01} and~\oldref{TableComparisonEx01} is also in the line with Theorem~\oldref{FalseDiscovery}.
\end{ex}

\begin{ex}\label{Ex02}
In this setting, we consider the DAG structure in Figure~\oldref{fig:Ex02} with the corresponding relationship between the nodes described in the following set of Equations~\ref{eq:Ex02}. Node $\bbx_5$ in this DAG does not satisfy the tree-neighborhood property, and therefore, DAG-FOCI is not able to recover $\{2, 3\}$ as the set of parents of $\bbx_5$. Nevertheless, from Theorem~\oldref{FalseDiscovery} we expect that with a large enough sample size, DAG-FOCI does not output any false-positives. To study this result empirically, for each sample size $n = 1000, 2000, \cdots, 8000$, we have generated 100 samples according to~\ref{eq:Ex02}. We have counted the number of times that DAG-FOCI has resulted in at least one false-positive detection. Table~\oldref{TableEx02} summarizes the results of these simulations. As it can be seen, as $n$ grows, the number of trials with at least one false-positive converges to zero. This observation is in line with Theorem~\oldref{FalseDiscovery} and shows the conservative nature of DAG-FOCI even when the tree-neighborhood assumption is not valid.
\begin{figure}[ht]
\centering
\begin{tikzpicture}
[->, node/.style={}]
\node[node]    (X1)                     {$\bbx_1$};
\node[node]    (X3)    [below= of X1]   {$\bbx_3$};
\node[node]    (X2)    [left= of X3]   {$\bbx_2$};
\node[node]    (X4)    [right= of X3]   {$\bbx_4$};
\node[node]    (X5)    [below= of X3]   {$\bbx_5$};
\node[node]    (X6)     [below = of X5]  {$\bbx_6$};
\node[node]    (X7)    [below right= of X5]   {$\bbx_7$};
\path (X1) edge (X2);
\path (X1) edge (X3);
\path (X1) edge (X4);
\path (X2) edge (X5);
\path (X3) edge (X5);
\path (X4) edge (X6);
\path (X5) edge (X6);
\path (X4) edge (X7);
\path (X5) edge (X7);
\end{tikzpicture}
\caption{DAG in Setting~\ref{Ex02}. The goal is to infer the parents of $\bbx_5$.}
\label{fig:Ex02}
\end{figure}
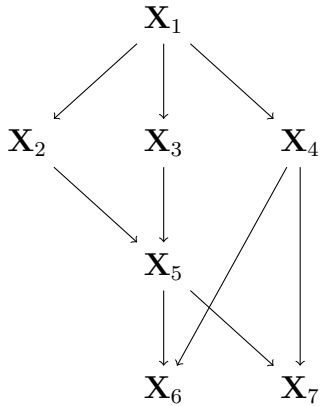
\begin{equation}\label{eq:Ex02}
\begin{split}
&\bbx_1, \varepsilon_i \overset{\text{i.i.d}}{\sim}\mathcal{N}(0, 1), \\
&\bbx_2 = \arctan(\bbx_1 + \varepsilon_2), \qquad\bbx_5 = \arctan(\bbx_2 + \bbx_3 + \varepsilon_5),\\
&\bbx_3 = \arctan(\bbx_1 + \varepsilon_3), \qquad\bbx_6 = \arctan(\bbx_4 + \bbx_5 + \varepsilon_6),\\
&\bbx_4 = \arctan(\bbx_1 + \varepsilon_4), \qquad\bbx_7 = \arctan(\bbx_4 + \bbx_5 + \varepsilon_7).
\end{split}
\end{equation}
\begin{table}[h!]
\begin{center}
\centering
\def\arraystretch{1}
\begin{tabular}{lrrrrrrrr}
  \hline
  $n$ & 1000 & 2000 & 3000 & 4000 & 5000 & 6000 & 7000 & 8000 \\
  \hline
  number of runs with false positive & 16 & 9 & 4 & 3 & 1 & 0 & 0 & 0 \\
  \hline
\end{tabular}
\end{center}
\caption{Number of runs that DAG-FOCI outputs at least one false discovery for node $\bbx_5$ in Setting~\ref{Ex02} for different sample sizes $n$, out of 100 independent simulation runs.}
\label{TableEx02}
\end{table}
\end{ex}

\subsection{Application to Real Data}\label{Real Data}
We consider the dataset in Sachs et al.~\cite{sachs2005causal}. This dataset consists of 7466 measurements of the abundance of phosphoproteins and phospholipids recorded under different experimental conditions in primary human immune system cells. Following the setting in \cite{wang2017permutation}, we consider only 5846 measurements in which the perturbations of receptor enzymes are identical. The observational distribution then is defined as the model where only the receptor enzymes are perturbed. This results in 1755 observational measurements and 4091 interventional measurements across five environments.  Table~\oldref{table:cytometry_summary} displays the number of samples for the observational environment (no interventions on any of the proteins) as well as for the five interventional environments; we also present the proteins that have been identified to receive an intervention (well approximated by a do-intervention) for each interventional environment \cite{wang2017permutation}.

\begin{table}[h!]
\def\arraystretch{1}
    \begin{center}
    \begin{tabular}{lllllll}
    \hline
    Intervention & None & Akt & PKC & PIP2 & Mek & PIP3 \\ \hline
    $\#$ samples & 1755 & 911 & 723 & 810 & 799 & 848  \\
    \hline
    \end{tabular}
    \end{center}
\caption{\small{Number of samples under each protein intervention for the flow cytometry dataset from Sachs et al.~\cite{sachs2005causal}}.}
\label{table:cytometry_summary}
\end{table}
{\bf Analysis with observational  environment}: We first apply DAG-FOCI (in  Algorithm~\oldref{alg:overall}) on observational data for the target variables $\{\text{Akt,PKC,PIP2, MEK, PIP3}\}$. The second column of Table~\oldref{table:summaryCyt} displays the output $\hat{\mathcal{P}}_n(Y)$ containing the estimated sets of parents for each target variable of interest. Examining these results, we can conclude the undirected edges $\text{PKC}-\text{MEK}$, $\text{PIP2}-\text{MEC}$, $\text{PIP3} -\text{ERK}$, $\text{MEK}-\text{PIP3}$, $\text{MEK}-\text{Raf}$, where an undirected edge between a pair of variables means that they are connected in the skeleton of the directed graph, but the direction of the edge is not identifiable by DAG-FOCI.  

{\bf Combining observational and interventional environments}:
We next exploit interventional environments to infer directionality in some of the edges found in our observational analysis. Specifically, we apply the extension of DAG-FOCI to interventional data (in Algorithm~\oldref{alg:interven}) to target variables $\{\text{Akt,PKC,PIP2, MEK, PIP3}\}$; as shown in Table~\oldref{table:cytometry_summary}, each of these variables has received a do-intervention and thus combining observational and interventional environments can potentially be useful for improving identifiability. The second column in Table~\oldref{table:summaryCyt} displays the estimated Markov boundary of each target variable in the corresponding interventional environment (step 2 in Algorithm~\oldref{alg:interven}) and the last two columns in Table~\oldref{table:summaryCyt} show the output of Algorithm~\oldref{alg:interven}: the sets $\tilde{\mathcal{P}}_n(Y),\tilde{C}_n(Y)$ consisting of all possible parental sets and child nodes of a target variable after incorporating knowledge of the interventions. We identify a child node $\tilde{C}_n(Y) = \{\text{Raf}\}$ corresponding to a target variable MEK and thus resolve the directionality between the nodes RaF and MEK. The remaining edges remain unresolved. A summary of all the causal conclusions made by DAG-FOCI are presented in Table~\oldref{table:EdgeCyt}.

\begin{table}[h!]
\centering
\tabcolsep=0.09cm
\def\arraystretch{1}
\begin{tabular}{r|c|c|c|c}
  \hline
 &  DAG-FOCI's Output $\hat{\mathcal{P}}_n(Y)$ & $\widehat{\MB}_{I}(Y)$ & $\tilde{\mathcal{P}}_n(Y)$ & $\tilde{C}_n(Y)$\\
 \hline
 Akt &  $\emptyset$ & \{Erk, PKA\} & $\emptyset$ & $\emptyset$\\
 PKC & \{\{Mek\}, $\emptyset$\} & \{P38\} &  \{\{Mek\}, $\emptyset$\} & $\emptyset$ \\
 PIP2 & \{\{Mek\}, $\emptyset$\}& \{PKC\} &  \{\{Mek\}, $\emptyset$\} & $\emptyset$ \\
 PIP3 & \{\{Erk\}, $\emptyset$\} & \{PIP2, Plc$\gamma$\} &  \{\{Erk\}, $\emptyset$\} & $\emptyset$\\
 Mek & \{\{PIP3\}, \{Raf\}, $\emptyset$\} & \{Raf, Erk\} &\{\{PIP3\}, $\emptyset$\} & \{Raf\} \\
  \hline
\end{tabular}
\caption{The first column presents the list of target variables we analyze. The second column contains the output of DAG-FOCI Algorithm~\oldref{alg:overall} $\hat{\mathcal{P}}_n(Y)$ on observational data. The third column displays the Markov boundary $\widehat{\MB}_I(Y)$ of a target variable in the interventional environment with do-intervention on the target variable. The last two columns are the outputs of modification of DAG-FOCI in Algorithm~\oldref{alg:interven} after combining observational and interventional data.}
\label{table:summaryCyt}
\end{table}

\begin{table}[h!]
\centering
\tabcolsep=0.09cm
\def\arraystretch{1}
\begin{tabular}{lrrrrrrr}
  \hline
  & \cite{taeb2021perturbations} & \cite{sachs2005causal} & \cite{mooij2013cyclic} & \cite{meinshausen2016methods} & \cite{eaton2007exact} & \cite{wang2017permutation} & \cite{triantafillou2015constraint}\\
 \hline
 Mek$\rightarrow$ Raf & \checkmark & x & \checkmark & \checkmark & \checkmark & x & \checkmark\\
 PKC - \text{Mek} &$\rightarrow$ &$\rightarrow$   & $\rightarrow$ & x & $\rightarrow$  & x & x\\
 {PIP2} - \text{Mek}& x & x & x & x & x & x & x\\
 \text{PIP3} - \text{Erk}& x & x & x & x & x & x & x\\
 \text{Mek} - \text{PIP3}& x & x & x & x & x & x & x\\
 \hline
\end{tabular}
\caption{Comparing the discovered edges (directed and undirected,i.e., edge in the skeleton) by using DAG-FOCI with observational and interventional data (left column) with different causal discovery methods (with indicated reference). The symbol ``$\checkmark$" means that the same directed edge is inferred, ``$\rightarrow$" indicates the edge direction inferred for an undirected edge from DAG-FOCI, and ``x" denotes that there is no overlap between (either directed or undirected) edges.}
\label{table:EdgeCyt}
\end{table}
\section{Discussion}
\label{sec:conclusion}
In this paper, we proposed DAG-FOCI, a computationally efficient algorithm to identify the set of parents of a target variable of interest. We provided finite-sample guarantees for DAG-FOCI and illustrated its utility in synthetic and real data experiments. A number of interesting directions for future investigation arise from our work. First, our theoretical analysis in Section~\oldref{sec:Theory} relies heavily on the `$\delta$-Markov Gap' assumption. It would be of interest to explore the extent to which this condition can be relaxed if we impose restrictions on the functional relationships and the error structures. Second, our aim in this paper was to infer accurate causal relations without making any distributional assumptions; the price we pay for this generality is that our method may produce few discoveries. An interesting research direction is to improve the power of DAG-FOC by, for example, developing a hybrid procedure that incorporates valid conditional independence relationships (when available) into DAG-FOCI. Finally, we assumed throughout that all relevant variables are observed. Extensions of DAG-FOCI that account for latent confounding would be of practical interest.  
\section*{Acknowledgments}
This project has received funding from the European Research Council (ERC) under the European Union’s Horizon 2020 research and innovation programme (grant agreement No. 786461).
\bibliography{ref}
\newpage
\appendix
\section{Proofs}\label{proofs}
\subsection{Proof of \ref{p1}-\ref{p4}}\label{Appendixp1p6}
\begin{proof} For any $k$ such that $|\texttt{pa}_k(Y)|\geq 2$ let $i, j\in\texttt{pa}_k(Y)$. We have $\bbx_i\not\perp\bbx_j\mid Y$, therefore $i\in\MB(\bbx_j)$. This shows the validity of~\ref{p1}. 

For $i\in\texttt{pa}_k(Y)$ and $j\in\MB(Y)\backslash\texttt{pa}_k(Y)$, if $Y$ satisfies the tree-neighborhood assumption, conditional on $Y$, $\bbx_i$ and $\bbx_j$ are independent. Also $Y$ belongs to the Markov boundary of $\MB(\bbx_i)$ and $\MB(\bbx_j)$, therefore $i\not\in\MB(\bbx_j)$. This proves~\ref{p3}.

For $i, j\in\MB(Y)$ such that $\bbx_i\in\MB(\bbx_j)$ and $\bbx_j\in\MB(\bbx_i)$ and $\bbx_i\perp\bbx_j$, we can only have $i, j\in\texttt{pa}_k(Y)$ for some $k\in[K]$. For the other direction note that by~\ref{p1} we have $\bbx_i\in\MB(\bbx_j)$ and $\bbx_j\in\MB(\bbx_i)$ and since $Y$ satisfies the tree-neighborhood assumption, the only path containing $i, j$ in the skeleton of $D^\star$ goes through $Y$. This implies $\bbx_i\perp\bbx_j$ which proves~\ref{p2}.

For $k\in[K]$ such that $|\texttt{pa}_k(Y)|>1$, if $Y$ satisfies the tree-neighborhood assumption, by \ref{p1} and \ref{p3} we cannot have $i, j\in[K]$ such that $\texttt{pa}_i(Y)\neq\texttt{pa}_j(Y)$ and therefore $K = 1$. On the other hand assume $K > 1$, and there exists some $k\in[K]$ such that $|\texttt{pa}_k(Y)| > 1$. Nodes $\texttt{pa}_k(Y)$ and $Y$ form a v-structure that should be present in all the Markov equivalence classes. This implies for all $i\in[K]$ we must have $\texttt{pa}_k(Y)\subseteq\texttt{pa}_i(Y)$. With the same argument we must have $\texttt{pa}_i(Y)\subseteq\texttt{pa}_k(Y)$ which means for all $i\in[K]$ we have $\texttt{pa}_i(Y)=\texttt{pa}_k(Y)$ which is in contradiction with $K > 1$. This concludes~\ref{p4}.
\end{proof}

\subsection{Proof of Theorem~\oldref{main}}
\begin{proof}
Under assumptions~\ref{a1}, and~\ref{a2}, Theorem 6.1 in~\cite{azadkia2021simple} gives us 
\[
\pp(\widehat{\MB}(Y)\text{ is a Markov blanket})\geq 1 - L_1 p^{L_2}e^{-L_3 n}.
\]
Hence we only need to show that with the additional assumption~\ref{a3}, $\widehat{\MB}(Y)$ is exactly the Markov boundary. 

Let $j_1, j_2, \ldots, j_p$ be the complete ordering of all variables produced by the step-wise algorithm in FOCI. Let $S_0 = \emptyset$, and $S_k = \{j_1, \ldots, j_k\}$ for all $k\leq p$. Let 
\begin{eqnarray*}
    Q(Y, \bbz) = \int\var(\pp(Y\geq t\mid \bbz))d\mu(t), 
\end{eqnarray*}
and 
\begin{eqnarray*}
    Q_n(Y, \bbz) = \frac{1}{n^2}\sum_{i=1}^n\big(\min\{R_i, R_{M(i)}\} - \frac{L_i^2}{n}\big),
\end{eqnarray*}
be the numerator of $T(Y, \bbz)$ and $T_n(Y, \bbz)$ respectively. Remember that $R_i = \sum_{j=1}^n1\{Y_j\leq Y_i\}$, and $L_i = \sum_{j=1}^n1\{Y_j\geq Y_i\}$, and $M(i)$ is the index of the nearest neighbor of $\bbz_i$ among all the other $\bbz_j$'s.

As defined in proof of Theorem 6.1 in~\cite{azadkia2021simple}, let $E^\prime$ be the event that for all $1\leq k\leq K$, $|Q_n(Y, \bbx_{S_k}) - Q(Y, \bbx_{S_k})|\leq\delta/8$ where $K$ is the integer part of $\/\delta + 2$. Note that $E^\prime$ and assumption~\ref{a3} together imply that members of $\MB(Y)$ appear first in the ordering by FOCI. This completes the proof.
\end{proof}

\subsection{Proof of Theorem~\oldref{bonferroni}}
\begin{proof}
We divide event $\mathcal{P}(Y)\not\subseteq\hat{\mathcal{P}}_n(Y)$, into two parts. 
\begin{eqnarray*}
     \pp(\mathcal{P}(Y)\not\subseteq\hat{\mathcal{P}}_n(Y)) &=& \pp(\{\mathcal{P}(Y)\not\subseteq\hat{\mathcal{P}}_n(Y)\}\cap \{\MB(Y) = \widehat{\MB}(Y)\}) + \\
     && \pp(\{\mathcal{P}(Y)\not\subseteq\hat{\mathcal{P}}_n(Y)\}\cap \{\MB(Y) \neq \widehat{\MB}(Y)\}) \\
     &\leq & \pp(\{\mathcal{P}(Y)\not\subseteq\hat{\mathcal{P}}_n(Y)\}\cap \{\MB(Y) = \widehat{\MB}(Y)\}) + \\
     && \pp(\{\MB(Y) \neq \widehat{\MB}(Y)\}).
\end{eqnarray*}
Theorem~\oldref{main} gives an upper bound on $\pp(\{\MB(Y) \neq \widehat{\MB}(Y)\})$. Therefore we need to bound $\pp(\mathcal{P}(Y)\not\subseteq\{\hat{\mathcal{P}}_n(Y)\}\cap \{\MB(Y) = \widehat{\MB}(Y)\})$. When $\MB(Y) = \widehat{\MB}(Y)$, then the estimation of parents can go wrong if at least one of the following happens
\begin{enumerate}[label=(F.\arabic*),ref=F.\arabic*]
    \item \label{F1} For at least one $i\in\MB(Y)$, $\widehat{\MB}(\bbx_i)\neq\MB(\bbx_i)$,
    \item \label{F2} We falsely reject the independence between at least two parental nodes,
    \item \label{F3} We falsely accept the independence between all the non-independent pairs in a non-parental component. 
\end{enumerate}
Theorem~\oldref{main} gives us a bound to control~\ref{F1}. Note that controlling~\ref{F2} is just controlling the Type I error in our independence test procedure.  

Assume that for all $i\in\MB(Y)$, we get $\widehat{\MB}(\bbx_i) = \MB(\bbx_i)$. In this case our estimated ${G}^\star_Y$ is equal to its population version and it has $R$ disjoint connected component. Consider a non-parental component $G_i$ in ${G}^\star_Y$, containing child node $j$ and spouse nodes $i_1, \cdots, i_k$. DAG-FOCI would falsely recognize $G_i$ as a parental component if and only if all the p-values of independence tests between pairs of nodes in $G_i$ are bigger than $\alpha$. This means that for all the non-independent pairs $(\bbx_j, \bbx_{i_\ell})$ the independence test fail to reject the null hypothesis. For node $Y$ with the tree-neighborhood property we expect $|V(G_i)| - 1$ non-independent pairs of nodes in $G_i$, all the $(\bbx_{j}, \bbx_{i_\ell})$'s where $j$ is the child of $Y$ and $i_\ell$ is a spouse. The probability of this event is $(\beta(n))^{|V(G_i)| - 1}$, where $\beta(n)$ is defined in~\ref{betan}. Putting all these together using union bound gives us 
\begin{eqnarray*}
    \pp(\mathcal{P}(Y)\not\subseteq\hat{\mathcal{P}}_n(Y)) &\leq & (|\MB(Y)|+1)(L_1 p^{L_2}e^{-L_3 n}) + \binom{\max|\texttt{pa}_k(Y)|}{2}\alpha +\\
    && \sum_{i=1}^R\beta(n)^{|V(G_i)| - 1}1\{|V(G_i)|\geq 2\}  \\
    &\leq & (|\MB(Y)|+1)(L_1 p^{L_2}e^{-L_3 n}) + \binom{\max|\texttt{pa}_k(Y)|}{2}\alpha + \\
    && (R - 1)\beta(n)^{\min_{|V(G_i)|\neq 1}|V(G_i)| - 1} \\
    &\leq & (|\MB(Y)|+1)(L_1 p^{L_2}e^{-L_3 n}) + \binom{\max|\texttt{pa}_k(Y)|}{2}\alpha + (R - 1)\beta(n)
\end{eqnarray*}
Since $R \leq |\MB(Y)|$ we have
\begin{eqnarray*}
    \pp(\mathcal{P}(Y)\not\subseteq\hat{\mathcal{P}}_n(Y))&\leq& (|\MB(Y)|+1)(L_1 p^{L_2}e^{-L_3 n}) + \alpha\max_k|\texttt{pa}_k(Y)|(\max_k|\texttt{pa}_k(Y)| - 1)/2 +\\
    && (|\MB(Y)| - 1)\beta(n).
\end{eqnarray*}
This completes the proof of Theorem~\oldref{bonferroni}.
\end{proof}

\subsection{Proof of Lemma~\oldref{prop1}}
\begin{proof} 
Let $S$ be a set such that $S\subsetneq\MB(Y)$. Let $z\not\in\MB(Y)$ be an arbitrary node. Take $w\in\MB(Y)$ the closest member of $\MB(Y)$ on the path between $Y$ and $z$ in $D^\star$ to $Y$. Since $Y$ has the tree-neighborhood property this path is unique and $w$ is well-defined. First let's consider the case where $w\not\in S$. If the path from/to $Y$ to/from $z$ is directed in $D^\star$ then given $\bbx_w$, $Y$ and $\bbx_z$ are independent of each other. Since $w\in\MB(Y)$, this implies
\[
\varepsilon_z := Q(Y, (\bbx_S,\bbx_w)) - Q(Y, (\bbx_S, \bbx_z)) > 0.
\]
If the path between $Y$ and $z$ in $D^\star$ is not directed let $w^\prime$ be the farthest node on the path between $w$ and $z$ to $w$ such that the path between $w$ and $w^\prime$ is directed. 
\begin{enumerate}[label=(m.\arabic*),ref=(m.\arabic*)]
    \item \label{m1} $Y\leftrightarrows w\rightarrow\cdots\rightarrow w^\prime\leftarrow\cdots z$
    \item \label{m2} $Y\leftrightarrows w\leftarrow\cdots\leftarrow w^\prime\rightarrow\cdots z$
\end{enumerate}
Let $v = \{w^\prime, \cdots, z\}$ the set consist of all the nodes on the path between $w^\prime$ and $z$. Then we replace this part of the path by the single node $v$. This will modify~\oldref{m1} and~\oldref{m2} into
\begin{enumerate}[label=(n.\arabic*),ref=(n.\arabic*)]
    \item \label{n1} $Y\leftrightarrows w\rightarrow\cdots\rightarrow v$
    \item \label{n2} $Y\leftrightarrows w\leftarrow\cdots\leftarrow v$
\end{enumerate}
When the path between $Y$ and $v$ is directed again we have
\[
Q(Y, (\bbx_S, \bbx_w)) - Q(Y, (\bbx_S, \bbx_v)) > 0.
\]
By monotonocity of $Q$ in its second argument we have $Q(Y, (\bbx_S, \bbx_v))\geq Q(Y, (\bbx_S, \bbx_z))$ and this results 
\[
\varepsilon_z := Q(Y, (\bbx_S, \bbx_w)) - Q(Y, (\bbx_S, \bbx_z)) > 0.
\]
When the path between $Y$ and $v$ is not directed but $w$ is the parent of $Y$, i.e., $Y\leftarrow w\rightarrow\cdots v$, we have $Y\perp \bbx_v\mid \bbx_w$ and similar to the previous case we get 
\[
\varepsilon_z := Q(Y, (\bbx_S, \bbx_w)) - Q(Y, (\bbx_S, \bbx_z)) > 0.
\]
So the only remaining case is $Y\rightarrow w\leftarrow\cdots v$. In this case there exists $w_2\in\MB(Y)$ such that $w\leftarrow w_2$ and hence $w_2$ is a spouse of $Y$. Then we have $Y\perp\bbx_v\mid\bbx_{w_2}$ and therefore 
\[
\varepsilon_z := Q(Y, (\bbx_S, \bbx_{w_2})) - Q(Y, (\bbx_S, \bbx_z)) > 0.
\]

If $w\in S$, then 
\[
Q(Y, \bbx_S) = Q(Y, (\bbx_S, \bbx_w)) = Q(Y, (\bbx_S, \bbx_z)). 
\]
Since $\bbx_S$ is a strict subset of $\MB(Y)$, for any $w^\prime\in\MB(Y)\setminus S$ we have
\[
Q(Y, (\bbx_S, \bbx_{w^\prime})) - Q(Y, (\bbx_S, \bbx_z)) > 0.
\]
Let 
\[
\varepsilon_z := \max_{{w^\prime}\in\MB(Y)\setminus S}Q(Y, (\bbx_S, \bbx_w)) - Q(Y, (\bbx_S, \bbx_z)).
\]
Therefore $\min(\delta, 4\underset{{z\not\in\MB(Y), S\subsetneq\MB(Y)}}{\min}\{\varepsilon_z\})$ satisfies the value for the $\delta$-Markov Gap property. This finishes proof of Lemma~\oldref{prop1}.
\end{proof}

\subsection{Proof of Theorem~\oldref{FalseDiscovery}}
\begin{proof}
Note that under assumptions~\ref{a1}-\ref{a3} with high probability we have $\widehat{\MB}(Y) = \MB(Y)$. Also under assumptions~\ref{a1}-\ref{a2} for nodes $z\in\MB(Y)$ with high probability we have $\MB(\bbx_z)\subseteq\widehat{\MB}(\bbx_z)$. Therefore to find an upper bound on $\pp(F)$ we can decompose this event into two disjoint part.
\begin{eqnarray*}
\pp(F) &=& \pp(F\cap\{\forall z\in\MB(Y): \MB(\bbx_z)\subseteq\widehat{\MB}(\bbx_z)\}\cap\{\widehat{\MB}(Y) = \MB(Y)\}) +\\
&& \pp(F\cap(\{\exists z\in\MB(Y): \MB(\bbx_z)\not\subseteq\widehat{\MB}(\bbx_z)\}\cup\{\widehat{\MB}(Y) \neq\MB(Y)\})).
\end{eqnarray*}
By Theorem~\oldref{main},~\cite[Th.6.1]{azadkia2021simple}, and union bound we have 
\begin{eqnarray*}
& &\pp(F\cap(\{\exists z\in\MB(Y):\MB(\bbx_z)\not\subseteq\widehat{\MB}(\bbx_z)\}\cup\{\widehat{\MB}(Y)\neq\MB(Y)\})) \\
& &\leq (|\MB(Y)| + 1)L_1p^{L_2}e^{L_3n}.
\end{eqnarray*}
Now assume for all $i\in\MB(Y)$ we have $\MB(\bbx_i)\subseteq\widehat{\MB}(\bbx_i)$ and $\widehat{\MB}(Y) = \MB(Y)$, and $F$ occurs. Therefore there exists $S\in\hat{\mathcal{P}}_n(Y)$ and $w\in S$ such that $w\not\in\texttt{pa}_k(Y)$ for any $k\in[K]$. Note that in this case if $\texttt{pa}_k(Y)\cap S\neq\emptyset$ then $\texttt{pa}_k(Y)\subsetneq S$, since $\widehat{\MB}(Y) = \MB(Y)$ and for all $i\in\MB(Y)$ we have $\MB(\bbx_i)\subseteq\widehat{\MB}(\bbx_i)$, unless we mistakenly reject the null assumption of independence for a pair of nodes in $\texttt{pa}_k(Y)$. Now consider the following two cases:
\begin{enumerate}
    \item $\texttt{pa}_k(Y)\subsetneq S$. In this case $w$ is a descendent of $v$ where $v\in\texttt{pa}_k(Y)$ and therefore $\bbx_w\not\perp\bbx_v$. 
    \item $\texttt{pa}_k(Y)\cap S =\emptyset$. In this case $S$ contains at least one child or spouse of $Y$ like $v$. Again in this case $\bbx_w\not\perp\bbx_v$. 
\end{enumerate}
Therefore $F$ happens if at least for one pair of dependent pairs we wrongly accept the null hypothesis of independence or for two independent pair of nodes we mistakenly reject the null assumption of independence. Therefore
\begin{eqnarray*}
& &\pp(F\cap\{\forall z\in\MB(Y)\mid \MB(\bbx_z)\subseteq\widehat{\MB}(\bbx_z)\}\cap\{\widehat{\MB}(Y) = \MB(Y)\}) \\
& &\leq \frac{1}{2}|\MB(Y)|^2(\beta(n) + \alpha).
\end{eqnarray*}
Therefore we have 
\begin{eqnarray*}
\pp(F) \leq \frac{1}{2}|\MB(Y)|^2(\beta(n) + \alpha) + (|\MB(Y)| + 1)L_1p^{L_2}e^{L_3n}.
\end{eqnarray*}
\end{proof}
\subsection{Proof of Corollary~\oldref{thm:interv}}
\begin{proof}
If $\tilde{C}_n(Y) \not\subseteq \texttt{ch}(Y)$ then there exists $i\in\tilde{C}_n(Y)\setminus\texttt{ch}(Y)$. This happens only if $i\in\widehat{\MB}_I(Y)$ and there exists $S\in\hat{\mathcal{P}}_n(Y)$ such that $i\in S$. If $\widehat{\MB}_I(Y) = \MB_I(Y)$ and $\widehat{\MB}(Y) = \MB(Y)$ and $\MB(\bbx_i)\subseteq\widehat{\MB}(\bbx_i)$ for all $i\in\MB(Y)$, this happens only if $i\in\texttt{sp}(Y)$. In this case $S\not\in\mathcal{P}(Y)$ and if $|S|>1$ we must have wrongly accepted the null hypothesis of independence between $\bbx_i$ and some other node in $S$. Note that if $|S| = 1$ then estimation of the Markov blanket of one of the members of $\MB(Y)$ must have gone wrong. Using the union bound and Theorem~\oldref{main} gives us
\begin{eqnarray*}
\mathbb{P}(\tilde{C}_n(Y) \not\subseteq \texttt{ch}(Y)) &\leq & (|\MB(Y)|+1)(L_1 p^{L_2}e^{-L_3 n_\text{obs}}) +\\
&& |\MB(Y)|^2\beta(n_\text{obs}) + L_1p^{L_2}e^{-L_3 n_I}.
\end{eqnarray*}
Note that with the extra assumptions that $Y$ satisfies the tree-neighborhood assumption, and that there are no do-interventions on any child node of $Y$, when $\widehat{\MB}_I(Y) = \MB_I(Y)$ and $\widehat{\MB}(Y) = \MB(Y)$ and $\MB(\bbx_i)\subseteq\widehat{\MB}(\bbx_i)$ for all $i\in\MB(Y)$, we have $\tilde{\mathcal{P}}_n(Y)=\{\texttt{pa}(Y)\}$. Therefore we get the following bound exactly as before. 
\begin{eqnarray*}
\mathbb{P}(\tilde{C}_n(Y) \not\subseteq \texttt{ch}(Y) \cup \tilde{\mathcal{P}}_n(Y)\neq\{\texttt{pa}(Y)\}) &\leq & (|\MB(Y)|+1)(L_1 p^{L_2}e^{-L_3 n_\text{obs}}) +\\
&& |\MB(Y)|^2\beta(n_\text{obs}) + L_1p^{L_2}e^{-L_3 n_I}.
\end{eqnarray*}
\end{proof}

\section{Violation of the $\delta$-Markov Gap Assumption}\label{sec:NoTree}
\begin{figure}[h]
\centering
\begin{tikzpicture}
[->, node/.style={}]
\node[node]    (X3)                 {$\bbx_3$};
\node[node]    (X1)[left=of X3]     {$\bbx_1$};
\node[node]    (X2)[right=of X3]    {$\bbx_2$};
\node[node]    (Y) [below=of X3]    {$Y$};
\path (X3) edge (X1);
\path (X3) edge (X2);
\path (X1) edge (Y);
\path (X2) edge (Y);
\end{tikzpicture}
\caption{Markov boundary of $\bby$ is $\{1, 2\}$.} 
\label{fig:Cycle}
\end{figure}
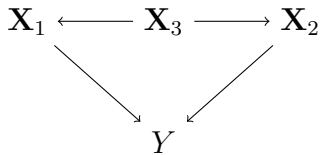

\begin{figure}[h!]
  \centering
  \includegraphics[width=0.8\textwidth]{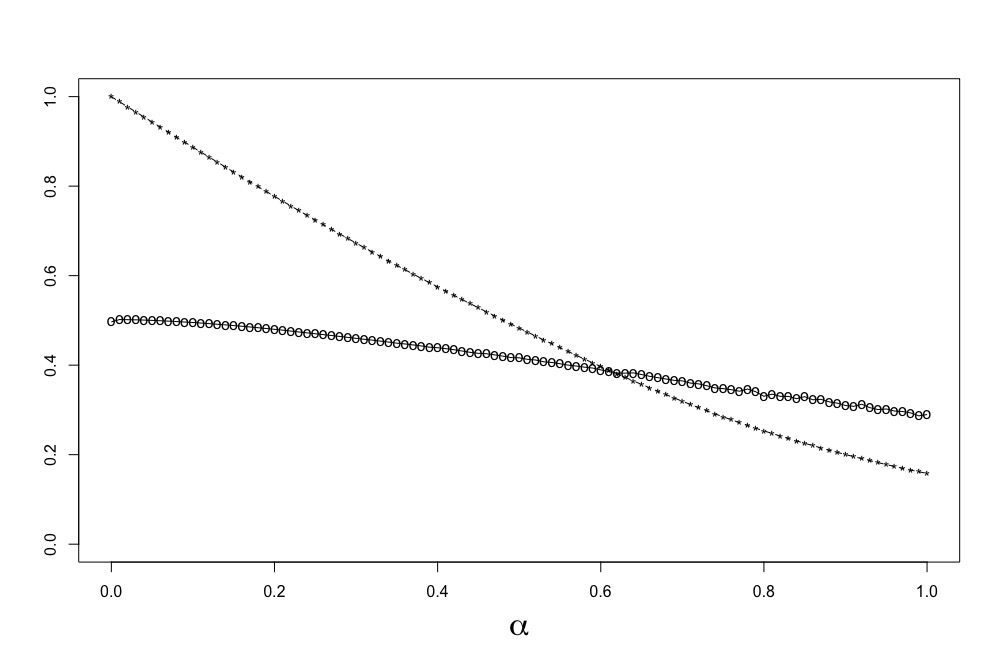}
  \caption{\small{The values of $T_n(Y, \bbx_3)$ ($*$) and $\max\{T_n(Y, \bbx_1), T_n(Y, \bbx_2)\}$ ($o$) for different values of $\alpha$, with $n = 10^4$. For small values of $\alpha$ we have $T_n(Y, \bbx_3) > \max\{T_n(Y, \bbx_1), T_n(Y, \bbx_2)\}$.}}
  \label{fig:alpha}
\end{figure}

Consider the DAG structure in Figure~\oldref{fig:Cycle}. The Markov boundary of $Y$ is $\{1, 2\}$. Note that since FOCI works in a greedy fashion, it is not hard to build distributions such that FOCI selects $3$ earlier than $1$ or $2$. For example let 
\begin{eqnarray}\label{violation1}
\begin{aligned}
&\bbx_1 = \sign(\bbx_3) + \alpha\varepsilon_1, \\
&\bbx_2 = |\bbx_3| + \alpha\varepsilon_2, \\
&Y = \bbx_1\bbx_2, 
\end{aligned}
\end{eqnarray}
where $\bbx_3, \varepsilon_1$, and $\varepsilon_2$ are i.i.d $\mathcal{N}(0, 1)$, and $\alpha$ is a non-negative constant. We assume for simplicity of exposition that $Y$ is noiseless. Figure~\oldref{fig:alpha} shows the values of $T_n(Y, \bbx_3)$ and $\max\{T_n(Y, \bbx_1), T_n(Y, \bbx_2)\}$ for $\alpha\in[0, 1]$, with $n = 10^4$, where $Y, \bbx_1, \bbx_2, \bbx_3$ follow the SEM in~\ref{violation1}. This shows that for small values of $\alpha$, the $\delta$-Markov Gap assumption~\ref{a3} is not satisfied and FOCI selects $3$ earlier than $1$, $2$. 
Although FOCI stops after selecting all $\{1, 2, 3\}$, therefore, returns a Markov blanket with high probability, this Markov blanket is no longer the Markov boundary. To get to the Markov boundary, one must perform further subset selection which makes the problem more computationally involved. 

\end{document}